\documentclass[a4paper,10pt,reqno]{amsart}
\usepackage{amssymb,latexsym,bbm,amsmath,amsfonts}
\usepackage{amsthm}
\usepackage[mathscr]{eucal}
\usepackage{rotating}
\usepackage{multirow}
\usepackage{slashbox}

\numberwithin{equation}{section}
\newtheorem{theorem}{Theorem}
\newtheorem{lemma}{Lemma}

\theoremstyle{definition}

\newtheorem{example}[theorem]{Example}

\theoremstyle{remark}

\begin{document}

\title[Growth rate delay dominated DDEs]
{Exact growth rates of solutions of delay--dominated differential equations with regularly varying coefficients}

\author{John A. D. Appleby}
\address{Edgeworth Centre for Financial Mathematics, School of Mathematical
Sciences, Dublin City University, Glasnevin, Dublin 9, Ireland}
\email{john.appleby@dcu.ie} \urladdr{webpages.dcu.ie/\textasciitilde
applebyj}

\author{Michael J. McCarthy}
\address{Edgeworth Centre for Financial Mathematics, School of Mathematical Sciences,
Dublin City University, Glasnevin, Dublin 9, Ireland.}
\email{michael.mccarthy29@mail.dcu.ie}

\author{Alexandra Rodkina}
\address{The University of the West Indies, Mona Campus, Department of Mathematics and Computer Science, Mona, Kingston 7, Jamaica.} \email{alexandra.rodkina@uwimona.edu.jm}

\thanks{The first author was partially funded by the Science Foundation Ireland
grant 07/MI/008 ``Edgeworth Centre for Financial Mathematics''.
Michael McCarthy was funded by The Embark Initiative operated by the Irish Research Council
for Science, Engineering and Technology (IRCSET) under the project ``Explosions of stochastic
delay differential equations in Finance".}
\date{18 July 2013}

\begin{abstract}
In this paper we determine the exact rate of growth of the solution of a deterministic delay differential equation in which the delayed term is regularly varying at infinity and dominates, and determine criteria to characterise this dominance.
The preservation of growth rates using a uniform step size Euler scheme is also discussed.
\end{abstract}

\maketitle

\section{Introduction}
This paper examines the growth rate of $x(t)\to\infty$ as $t\to\infty$ of solutions
of the delay differential equation
\begin{equation} \label{eq.intro}
x'(t)=f(x(t))+g(x(t-\tau)), \quad t>0,
\end{equation}
We establish criteria on the size of $g$ relative to $f$ under which the solution of the delay equation does \emph{not} grow like the solution of the ordinary differential equation $y'(t)=f(y(t))$. In broad terms, we focus on the cases when $g$ grows polynomially, and $f$ grows sublinearly, though the general theory extends to cover more rapidly growing $g$ as well, and even recovers the exact exponential growth and characteristic equation in the linear case.

In \cite{AppMcCartRod10}, we established \emph{general} results for the exact rate of growth of solutions of \eqref{eq.intro} in which the delay term in some sense asymptotically dominates the instantaneous term.
The general theorems are obtained by employing a \emph{constructive comparison principle} (see Appleby~\cite{App05} and Appleby and Buckwar~\cite{AppBuck09}, for example). The asymptotic results are restated here in Section 2.2, and their hypotheses explained. In these general theorems, the sufficient conditions which describe this dominance, as well as the rate of growth of solutions, depend on the existence of an auxiliary function $\phi$ obeying certain asymptotic properties. Apart from some examples, we do not attempt systematically in \cite{AppMcCartRod10} to  demonstrate that such an auxiliary function $\phi$ exists, nor did we indicate how it might be constructed.

In this paper, we show when $g$ is regularly varying at infinity with positive index, and $f$ is sufficiently small,
that the rate of growth of solutions of \eqref{eq.intro} can be determined in the form
\[
\lim_{t\to\infty} \frac{G(x(t))}{t}=\lambda>0,
\]
for some function $G$ that is known in terms of $g$. This is achieved because, for such classes of problems, the auxiliary function $\phi$ can be found, and the exact asymptotic behaviour determined by applying general results. In addition, we show that for an explicit Euler scheme with uniform step size $h>0$ that the asymptotic behaviour is
preserved, in the sense that for every $h$ there exist $0<\lambda(h)<+\infty$ such that
\[
\lim_{t\to\infty} \frac{G(x_h(t))}{t}=\lambda(h), %\quad \limsup_{t\to\infty} \frac{G(x_h(t))}{t}=\lambda_U(h)
\]
and $\lim_{h\to 0}\lambda(h)=\lambda$, % \lim_{h\to 0} \lambda_U(h)=\lambda$,
where $x_h$ is the extension to continuous time of the Euler scheme.

Statements and discussion of the main results in continuous--time, as well as examples, are given in Section \ref{sec.continuousgrowth}. The preservation of the asymptotic rate of growth is considered in Section \ref{sec.discretegrowth}. Proofs of continuous time results are deferred to Section \ref{sec.continuousproofs},
with results for the discretisations being supplied in Section~\ref{sec.discreteproofs}. We do not address here the asymptotic behaviour when the instantaneous term dominates.

\section{Statement and Discussion of Main Continuous Time--Results}  \label{sec.continuousgrowth}

\subsection{Notation and preliminary results: existence and non--explosion}
In this paper, $\mathbb{R}$ stands for the real numbers, $\mathbb{N}$ for the natural numbers and $\mathbb{Z}$ for the integers. A function $k:[0,\infty)\to (0,\infty)$ is said to be \emph{regularly varying at infinity} with index $\alpha\in\mathbb{R}$ if
$\lim_{x\to\infty} k(\lambda x)/k(x)=\lambda^\alpha$ for each $\lambda>0$.
We write $k\in \mbox{RV}_\infty(\alpha)$. The reader is referred to Bingham, Goldie and Teugels~\cite{BingGoldTeu89} for results on regularly varying functions. If $I$ and $J$ are subintervals of $\mathbb{R}$, the space $C(I;J)$ contains all continuous functions $\phi:I\to J$.

We make some hypotheses regarding our problem. Suppose
\begin{equation} \label{eq:f1}
f\in C((0,\infty);[0,\infty)) \mbox{ is locally Lipschitz
continuous}.
\end{equation}
and obeys
\begin{equation} \label{eq:f2}
\int_{1}^\infty \frac{1}{f(u)}\,du=+\infty.
\end{equation}
We interpret this condition as being satisfied if $f$ is identically zero.
Suppose also that
\begin{equation} \label{eq:g}
g\in  C((0,\infty);(0,\infty)).
\end{equation}
Let $\tau>0$ and suppose that
\begin{equation} \label{eq:psi}
\psi\in C([-\tau,0];(0,\infty)),
\end{equation}
and consider the delay--differential equation given by
\begin{equation} \label{eq:cnseq}
x'(t)=f(x(t))+g(x(t-\tau)), \quad t>0; \quad x(t)=\psi(t), \quad t\in[-\tau,0].
\end{equation}
The following result then holds.
\begin{theorem} \label{thm:thmcnsnonexplosion}
Let $f$ obey (\ref{eq:f1}), (\ref{eq:f2}), $g$ obeys
(\ref{eq:g}), and $\psi\in C([-\tau,0];(0,\infty))$ where
$\tau>0$. Then there is $x\in C([-\tau,\infty)$ which is the
unique continuous solution of (\ref{eq:cnseq}) and which moreover
obeys
\begin{equation} \label{eq:cnsnonexplode}
 \lim_{t\to \infty} x(t)= \infty.
\end{equation}
\end{theorem}
The condition (\ref{eq:f2}) prevents a finite time explosion. Note that $g$ being positive forces $x$ to be increasing on $[0,\infty)$, and this ensures that \eqref{eq:cnsnonexplode} holds, because $\lim_{t\to\infty} x(t)=:L\in (0,\infty)$ forces $g(L)=0$, a contradiction.

\subsection{Statement and Discussion of General Comparison Results}
Before we state our general comparison results, we first introduce some notation and auxiliary functions.
Since $\psi(t)>0$ for all $t\in [-\tau,0]$ we may define $\psi^\ast:=\max_{-\tau\leq s\leq 0}\psi(s)>0$. Suppose that $\phi:(\psi^\ast,\infty)\to(0,\infty)$ is continuous, and define
\begin{equation} \label{def.Gamma}
\Gamma(x)=\int_{\psi^\ast}^x \frac{1}{\phi(u)}\,du, \quad x>\psi^\ast.
\end{equation}
Suppose that
\begin{equation} \label{eq.Gammatoinfty}
\lim_{x\to\infty}\Gamma(x)=+\infty.
\end{equation}
Define also for $c>0$ the function $\Gamma_c$ given by
\begin{equation} \label{def.Gammac}
\Gamma_c(x) = \frac{1}{c}\Gamma(x), \quad x>\psi^\ast.
\end{equation}
In our first main result, which appears as Theorem 1 in \cite{AppMcCartRod10}, we claim that if the delayed term $f$ is asymptotically dominated by the instantaneous term $g$, then the solution of (\ref{eq:cnseq}) behaves according to the ordinary differential equation  $z'(t)=\phi(z(t))$.
\begin{theorem} \label{thm:limsup}
Suppose that $f$ obeys (\ref{eq:f1}) and (\ref{eq:f2}).
Let $g$ be non--decreasing and obey (\ref{eq:g}) and let $\tau>0$
and $\psi\in C([-\tau,0];(0,\infty))$. Suppose that there exists a continuous function $\phi$ such that $\Gamma$, $\Gamma_c$ are defined by (\ref{def.Gamma}) and (\ref{def.Gammac}) respectively, and that $\Gamma$ obeys (\ref{eq.Gammatoinfty}). Suppose also that
\begin{equation} \label{eq.etaepstoeta}
\lim_{\epsilon\to 0^+}\eta(\varepsilon)=\eta,
\end{equation}
and suppose that
\begin{eqnarray} \label{eq.fsmallg}
\lim_{x\to\infty} \frac{f(x)}{g(\Gamma_{\eta(\varepsilon)}^{-1}(\Gamma_{\eta(\varepsilon)}(x)-\tau))}&=&0, \quad \mbox{for every $\varepsilon\in (0,1)$},\\
\label{eq.gphigamma}
\limsup_{x\to\infty} \frac{g(x)}{\phi(\Gamma_{\eta(\varepsilon)}^{-1}(\Gamma_{\eta(\varepsilon)}(x)+\tau))}&=&\bar{\eta}_\varepsilon\in [0,\infty)
\quad \mbox{for every $\varepsilon\in (0,1)$},
\end{eqnarray}
where
\begin{equation}\label{eq.naretaeta}
\sup_{\epsilon\in (0,1)} \bar{\eta}_\epsilon=:\bar{\eta}<\eta.
\end{equation}
If $x$ is the unique continuous solution of (\ref{eq:cnseq}), then
\begin{equation} \label{eq.growthlimsup}
\limsup_{t\to\infty} \frac{\Gamma(x(t))}{t}\leq \eta.
\end{equation}
\end{theorem}
We comment briefly on Theorem~\ref{thm:limsup} and its hypotheses. First, we note that the existence of a function $\phi$ obeying (\ref{eq.fsmallg}) and (\ref{eq.gphigamma}) is not assured by the theorem; the existence or construction of such a function must be achieved independently. However, it can be seen that (\ref{eq.gphigamma}) describes an asymptotic relationship between $\phi$ and $g$ only, and this is what identifies candidates for $\phi$. In the next section, we show that for a wide class of $g$ that suitable $\phi$ can be chosen. The condition (\ref{eq.fsmallg}) characterises the fact that the instantaneous term $f$ is dominated by the delayed term.

We now offer an improvement on Theorem~\ref{thm:limsup}. In it the condition \eqref{eq.fsmallg} is relaxed. In later examples we show that this enables asymptotic estimates to be extended to a wider class of problems.  
\begin{theorem} \label{thm:limsup2}
Suppose that $f$ obeys (\ref{eq:f1}) and (\ref{eq:f2}).
Let $g$ be non--decreasing and obey (\ref{eq:g}) and let $\tau>0$
and $\psi\in C([-\tau,0];(0,\infty))$. Suppose that there exists a continuous function $\phi$ such that $\Gamma$, $\Gamma_c$ are defined by (\ref{def.Gamma}) and (\ref{def.Gammac}) respectively, and that $\Gamma$ obeys (\ref{eq.Gammatoinfty}). Suppose also that \eqref{eq.etaepstoeta}
%\begin{equation} \label{eq.etaepstoeta}
%\lim_{\epsilon\to 0^+}\eta(\varepsilon)=\eta,
%\end{equation}
and suppose that $f$ obeys
\begin{equation} \label{eq.fsmallg2}
\lim_{x\to\infty} \frac{f(x)}{\phi(x)}=0,
\end{equation}
and that $g$ and $\phi$ obey \eqref{eq.gphigamma} where $\bar{\eta}_\varepsilon$ obeys \eqref{eq.naretaeta}
If $x$ is the unique continuous solution of (\ref{eq:cnseq}), then it obeys \eqref{eq.growthlimsup}. 
\end{theorem}

We now state a corresponding result which enables us to determine a lower bound on the rate of growth of solutions.
It appeared as Theorem 2 in~\cite{AppMcCartRod10}.
\begin{theorem} \label{thm:liminf}
Suppose that $f$ obeys (\ref{eq:f1}) and (\ref{eq:f2}).
Let $g$ be non--decreasing and obey (\ref{eq:g}) and let $\tau>0$
and $\psi\in C([-\tau,0];(0,\infty))$. Suppose that there exists a continuous function $\phi$ such that $\Gamma$, $\Gamma_c$ are defined by (\ref{def.Gamma}) and (\ref{def.Gammac}) respectively, and $\Gamma$ obeys (\ref{eq.Gammatoinfty}). Suppose also that
\begin{equation} \label{eq.muepstomu}
\lim_{\epsilon\to 0^+}\mu(\varepsilon)=\mu,
\end{equation}
and that $g$ and $\phi$ obey
\begin{equation}  \label{eq.gphigamma2}
\liminf_{x\to\infty} \frac{g(x)}{\phi(\Gamma_{\mu(\varepsilon)}^{-1}(\Gamma_{\mu(\varepsilon)}(x)+\tau(1-\epsilon)))}
=\bar{\mu}_\varepsilon\in (0,\infty] \quad \mbox{for every $\varepsilon\in (0,1)$},
\end{equation}
where
\begin{equation}\label{eq.narmumu}
\inf_{\epsilon\in (0,1)} \bar{\mu}_\epsilon=:\bar{\mu}>\mu.
\end{equation}
If $x$ is the unique continuous solution of (\ref{eq:cnseq}), then
\begin{equation} \label{eq.growthliminf}
\liminf_{t\to\infty} \frac{\Gamma(x(t))}{t}\geq \mu.
\end{equation}
\end{theorem}
As in Theorem~\ref{thm:limsup}, in which the condition (\ref{eq.gphigamma}) determines a relationship between $\phi$ and $g$, in Theorem~\ref{thm:liminf} there is a corresponding and closely related condition (\ref{eq.gphigamma2}) which describes the relationship between $g$ and $\phi$.

Contingent on other hypotheses being satisfied, we notice that the lower bound (\ref{eq.growthliminf}) and the upper bound (\ref{eq.growthlimsup}) incorporate the same function $\Gamma$. Therefore, under certain conditions we may combine Theorems~\ref{thm:limsup} and \ref{thm:liminf} to arrive at the exact asymptotic behaviour of $x$. This is the subject of the next result, which improves on a result in~\cite{AppMcCartRod10}.
\begin{theorem} \label{thm:lim}
Suppose $f$ obeys (\ref{eq:f1}) and (\ref{eq:f2}).
Let $g$ be non--decreasing and obey (\ref{eq:g}) and let $\tau>0$
and $\psi\in C([-\tau,0];(0,\infty))$. Suppose that there exists a continuous function $\phi$ such that $\Gamma$, $\Gamma_c$ are defined by (\ref{def.Gamma}) and (\ref{def.Gammac}), and that $\Gamma$ obeys (\ref{eq.Gammatoinfty}). Suppose also that there is $\eta>0$ such that $\mu(\epsilon)\to\eta$ and $\eta(\epsilon)\to \eta$ as $\epsilon\to 0$  and that $f$, $g$, and $\phi$ obey
(\ref{eq.fsmallg2}), (\ref{eq.gphigamma}) and (\ref{eq.gphigamma2}), where
\begin{equation}\label{eq.barmueta}
\sup_{\epsilon\in (0,1)} \bar{\eta}_\epsilon=:\bar{\eta}<\eta, \quad  \inf_{\epsilon\in (0,1)} \bar{\mu}_\epsilon=:\bar{\mu}>\eta.
\end{equation}
If $x$ is the unique continuous solution of (\ref{eq:cnseq}), then
\begin{equation} \label{eq.growthlim}
\lim_{t\to\infty} \frac{\Gamma(x(t))}{t}= \eta.
\end{equation}
\end{theorem}
Provided that a function $\phi$ can be found so that all the hypotheses of Theorem~\ref{thm:lim} are satisfied, the conclusion of Theorem~\ref{thm:lim} (viz., (\ref{eq.growthlim})) which describes an exact rate of growth, is sharp. 
%However, it is not certain whether the condition (\ref{eq.fsmallg}), which places a limit on the growth rate of $f$ relative to $g$, is sharp. However, it is clearly satisfied if $f$ grows sufficiently slowly. Therefore, in the next Section, in which concrete applications of Theorems~\ref{thm:limsup}--\ref{thm:lim} are presented, we have made conservative and explicit assumptions on the rate of growth of $f$.

\subsection{Application to equations with regularly varying $g$}
We consider some cases in which the unknown auxiliary function $\phi$ (and therefore $\Gamma$) in Theorems~\ref{thm:limsup}--\ref{thm:lim} can be constructed explicitly in terms of $g$. Essentially, our examples cover the cases where $g$ grows polynomially at either a sublinear or superlinear rate. First we consider the case where $g$ is in $\mbox{RV}_\infty(\beta)$ for $\beta\leq 1$ and $g(x)/x\to 0$ as $x\to\infty$.
\begin{theorem} \label{thm:rvbetalt1}
Let $f$ obey (\ref{eq:f1}), (\ref{eq:f2}). Let $g$ obey (\ref{eq:g}) be non--decreasing and let $\tau>0$
and $\psi\in C([-\tau,0];(0,\infty))$. Suppose $g\in \mbox{RV}_\infty(\beta)$ for some $\beta\leq 1$, $\lim_{x\to\infty}g(x)/x=0$, and $\lim_{x\to\infty}f(x)/g(x)=0$.
If $x$ is the unique continuous solution of (\ref{eq:cnseq}), then
\begin{equation} \label{eq.growthrvbetalt1}
\lim_{t\to\infty} \frac{1}{t}\int_1^{x(t)} \frac{1}{g(u)}\,du=1.
\end{equation}
\end{theorem}
This result is proven using Theorems~\ref{thm:limsup} and \ref{thm:liminf}; it recovers part (ii)
of Theorem 2.2 in Appleby, McCarthy and Rodkina~\cite{AppMcCartRod09}. Next we consider the case where $g$ is in $\mbox{RV}_\infty(1)$ but in which $g(x)/x\to \infty$ as $x\to\infty$, and use
Theorem~\ref{thm:lim} to determine the growth rate.
\begin{theorem} \label{thm:limrv1}
Let $f$ obey (\ref{eq:f1}), (\ref{eq:f2}). Let $g$ obey (\ref{eq:g}) and be non--decreasing. Let $\tau>0$ and $\psi\in C([-\tau,0];(0,\infty))$. Suppose $g\in \mbox{RV}_\infty(1)$, $x\mapsto g(x)/x$ is asymptotic to a non--decreasing function, $\lim_{x\to\infty}g(x)/x=\infty$, and 
\begin{equation} \label{eq.fsmallvsgrv1}
\lim_{x\to\infty} \frac{f(x)}{x\log(g(x)/x)}=0.
\end{equation}
%that there exists $\gamma<1$ such that
%\[
%\limsup_{x\to\infty} \frac{\log f(x)}{\log x}\leq\gamma<1.
%\]
Define
\begin{equation} \label{def:G}
G(x)=\int_{1}^x \frac{1}{u\log(1+g(u)/u)}\,du, \quad x>1.
\end{equation}
Then the unique continuous solution $x$ of (\ref{eq:cnseq}) obeys
\begin{equation} \label{eq.growthlimrv1}
\lim_{t\to\infty} \frac{G(x(t))}{t}=\frac{1}{\tau}.
\end{equation}
\end{theorem}
With a slightly stronger hypothesis on $f$ we can obtain the same conclusion on the growth rate, but by an alternative proof.
\begin{theorem} \label{thm:limnoinst}
Let $f$ obey (\ref{eq:f1}), (\ref{eq:f2}). Let $g$ obey (\ref{eq:g}) and be non--decreasing. Let $\tau>0$ and $\psi\in C([-\tau,0];(0,\infty))$. Suppose $g\in \mbox{RV}_\infty(1)$, $x\mapsto g(x)/x$ is asymptotic to a non--decreasing function, $\lim_{x\to\infty}g(x)/x=\infty$, and  $\lim_{x\to\infty}f(x)/x=0$.
If $G$ is defined by (\ref{def:G}), then the unique continuous solution $x$ of (\ref{eq:cnseq}) obeys
\begin{equation} \label{eq.growthlimdelay}
%\frac{1}{2\tau}\leq
\lim_{t\to\infty} \frac{G(x(t))}{t}=\frac{1}{\tau}.
%\limsup_{t\to\infty} \frac{G(x(t))}{t}\leq  \frac{1}{\tau}.
\end{equation}
\end{theorem}
The case where $g$ grows according to $g\in \mbox{RV}_\infty(\beta)$ for some $\beta\leq 1$ with $g(x)/x$ tending to a zero limit is covered by Theorem~\ref{thm:rvbetalt1}.

The proof of Theorem~\ref{thm:limnoinst} is facilitated by the following Lemma, which appears as Lemma 2.7 in Appleby, McCarthy and Rodkina~\cite{AppMcCartRod09}. It also motivates the choice of $\phi$ in Theorem~\ref{thm:limrv1}.
\begin{lemma} \label{lemma.l1}
Let $h>0$. Suppose $g\in C((0,\infty),(0,\infty))$, $g\in\mbox{RV}_\infty(1)$, $g(y)/y\to\infty$ as $y\to\infty$, and there is a function $g_1$ with $g_1(y)/g(y)\to1$ as $y\to\infty$ such that
$y\mapsto g_1(y)/y$ is non--decreasing. If $y_{n+1}=y_n+h g(y_n),\,n\geq 0$ and
$y_0=\xi>0$, then $\lim_{n\to\infty} G(y_n)/n=1$, where $G$ is defined by (\ref{def:G}).
\end{lemma}

If $g(x)/x$ tends to a finite non--zero limit, we are in the standard linear case, but even this is recovered independently of the standard linear theory by applying Theorems~\ref{thm:limsup} and \ref{thm:liminf}.
\begin{theorem} \label{ex:4}
Let $C>0$, $\tau>0$ and suppose that $\psi\in C([-\tau,0];(0,\infty))$.
Let $x$ be the unique continuous solution of (\ref{eq:cnseq}) with $f(x)/x\to 0$ and $g(x)/x\to C$ as $x\to\infty$. Then there is a unique $\lambda>0$ such that $\lambda=Ce^{-\lambda\tau}$ and $x$ obeys
$\lim_{t\to\infty} \log x(t)/t=\lambda$.
\end{theorem}

In the case when $g$ has a power--like growth faster which is faster than linear,
the rate of growth can be determined by means of Theorem~\ref{thm:lim}.
\begin{theorem} \label{thm:polyg}
Suppose that $f$ obeys (\ref{eq:f1}) and (\ref{eq:f2}). Let $g$ obey (\ref{eq:g}) be non--decreasing and let $\tau>0$ and $\psi\in C([-\tau,0];(0,\infty))$.
Suppose also that there exists $\beta>1$ such that $\lim_{x\to\infty} \log g(x)/\log x=\beta$ and 
\[
\lim_{x\to\infty} \frac{f(x)}{x\log x}=0.
\]
%$\limsup_{x\to\infty} \log f(x)/\log x\leq \gamma$. 
Then the unique continuous solution $x$ of (\ref{eq:cnseq}) obeys
\begin{equation} \label{eq.polygrowth}
\lim_{t\to\infty} \frac{\log\log x(t)}{t}=  \frac{\log (\beta)}{\tau}.
\end{equation}
\end{theorem}
The proofs of all these results are postponed to Section~\ref{sec.continuousproofs}.

\subsection{Examples}
We consider representatives example to which Theorem~\ref{thm:lim} can be applied.
%; these also appear in~\cite{AppMcCartRod10}.
%ref{thm:limnoinst} could be applied. However, we show that the estimates in (\ref{eq.growthlimdelay}) can sometimes be strengthened and replaced by a limit, which is secured by applying Theorem~\ref{thm:lim}.
%Our choice of $\phi$ in Theorem~\ref{thm:lim} is motivated by the form of $G$ in (\ref{def:G}), which plays the role of $\Gamma$ in Theorem~\ref{thm:limnoinst}.
For simplicity, we set $f$ to be identically zero.
%\begin{example} \label{ex:1}
%Suppose $g$ obeys (\ref{eq:g}) and is non--decreasing, and there exist $\alpha>0$ and $C>0$ such that $\lim_{x\to\infty} g(x)/(x(\log x)^\alpha)=C$, and $f(x)=0$ for all $x\geq 0$. Suppose $\tau>0$ and $\psi$ obeys (\ref{eq:psi}). Then the unique continuous solution $x$ of (\ref{eq:cnseq}) obeys $\lim_{t\to\infty} \log x(t)/(t\log t)=1/\tau$.
%\end{example}
%To see this, we note that $g$ obeys all the properties of Theorem~\ref{thm:limrv1}. Define $\phi(x)=x (\log_2 x)^2/(\log_2 x-1)$ for $x>e^e$. Then $\Gamma(x)=\log x/\log_2 x-e$ and
%$\lim_{x\to\infty} \log \Gamma^{-1}(x)/(x\log x)=1$. By Theorem~\ref{thm:limrv1} we have   $\lim_{t\to\infty}\Gamma(x(t))/t=1/\tau$. Thus $\lim_{t\to\infty}\log \Gamma^{-1}(\Gamma(x(t)))/(\Gamma(x(t))\log \Gamma(x(t)))=1/\tau$, which implies that $\lim_{t\to\infty}\log x(t)/(t\log t)=1/\tau$.
\begin{example}\label{ex:2}
Suppose $g$ obeys (\ref{eq:g}) and is non--decreasing, and there exists $C_1>0$ and $\alpha\in (0,1)$ such that
$\lim_{x\to\infty} g(x)/(x\exp((\log x)^\alpha))=C_1$, and $f(x)=0$ for all $x\geq 0$. Suppose $\tau>0$ and $\psi$ obeys (\ref{eq:psi}). Then the unique continuous solution $x$ of (\ref{eq:cnseq}) obeys
$\lim_{t\to\infty} \log x(t)/t^{1/(1-\alpha)}=(\eta(1-\alpha)/\tau)^{1/(1-\alpha)}$.
\end{example}
To see this, we note that $g$ obeys all the properties of Theorem~\ref{thm:limrv1}. For $x>e$ let $\phi(x)=x(\log x)^\alpha$. Then $\Gamma(x)=(\log(x)^{1-\alpha}-1)/(1-\alpha)$. By Theorem~\ref{thm:limrv1} we have   $\lim_{t\to\infty}\Gamma(x(t))/t=1/\tau$, which rearranges to give
$\lim_{t\to\infty} \log x(t)/t^{1/(1-\alpha)}=(\eta(1-\alpha)/\tau)^{1/(1-\alpha)}$.
%, \quad \log \Gamma^{-1}(x)=(1+(1-\alpha)\eta x)^{1/(1-\alpha)}.\]
%The calculations justifying these claims are postponed to the end of the next section.

We remark that the results can be applied to equations in which $g$ grows more rapidly than a polynomial function; here again is a representative example, which was considered without supporting calculations in~\cite{AppMcCartRod10}.
\begin{example}\label{ex:3}
Suppose $g$ obeys (\ref{eq:g}) and is non--decreasing, and there exists $C_1>0$ and $\alpha>1$ such that
$\lim_{x\to\infty} g(x)/\exp((\log x)^\alpha)=C_1$, and $f(x)=0$ for all $x\geq 0$. Suppose $\tau>0$ and $\psi$ obeys (\ref{eq:psi}). Then the unique continuous solution $x$ of (\ref{eq:cnseq}) obeys
$\lim_{t\to\infty} \log_3 x(t)/t=\log \alpha/\tau$.
\end{example}
To justify Example~\ref{ex:3}, set $\phi(x)=(1+x)\log(1+x)\log_2(1+x)$ for $x>e^e$. With $c:=\log_3(1+e^e)$, we have $\Gamma_\eta(x)=(\log_3(1+x)-c)/\eta$ and with $\lambda=e^{\eta\theta}$,
$\Gamma_\eta^{-1}(\Gamma_\eta(x)+\theta)=\exp((\log(1+x))^\lambda)-1$.
Therefore we have $\lim_{x\to\infty}
\phi(\Gamma_\eta^{-1}(\Gamma_\eta(x)+\theta))/(\exp([\log(1+x)]^\lambda)[\log x]^\lambda \log_2 x)=\lambda$.
Define $\eta(\epsilon)=(1+\epsilon)\log \alpha/\tau$ and $\mu(\epsilon)=\log \alpha/(\tau(1-\epsilon)^2)$. Then
\[
\lim_{x\to\infty} \frac{g(x)}{\phi(\Gamma_{\eta(\epsilon)}^{-1}(\Gamma_{\eta(\epsilon)}(x)+\tau))}=0
\]
and
\[
\lim_{x\to\infty} \frac{g(x)}{\phi(\Gamma_{\mu(\epsilon)}^{-1}(\Gamma_{\mu(\epsilon)}(x)+\tau(1-\epsilon)))}=\infty.
\]
Since $\eta(\epsilon), \mu(\epsilon)\to \log \alpha/\tau$ as $\epsilon\to\infty$, from Theorem~\ref{thm:lim} we have $\lim_{t\to\infty} \Gamma(x(t))/t=\log\alpha/\tau$, from which the result follows.

\section{Preservation of Growth Rates under Discretisation} \label{sec.discretegrowth}
Let $N\in \mathbb{N}$, and suppose that $h=\tau/N$. Consider the discretisation of \eqref{eq:cnseq}
according to
\begin{subequations} \label{eq:disceq}
\begin{align} \label{eq.disceqeq}
x_h(n+1)&=x_h(n)+hf(x_h(n))+hg(x_h(n-N)), \quad n\geq 0; \\
\label{eq.disceqic}
x_h(n)&=\psi(nh), \quad n=-N,\ldots,0.
\end{align}
\end{subequations}
We also find it of interest to define a continuous time extension of $x_h$.
If $(x_n)$ obeys \eqref{eq:disceq}, define $\bar{x}_h\in
C([-\tau,\infty),(0,\infty))$ by $\bar{x}_h(t)=\psi(t),\,t\in[-\tau,0]$,
\begin{equation} \label{def.lininterpolant}
\bar{x}_h(t)=x_n+(x_{n+1}-x_n)(t-nh)/h, \quad t\in[nh,(n+1)h], \quad
n\geq 0,
\end{equation}
so $\bar{x}_h$ takes the value $x_{n}(h)$ at time $nh$ for $n\geq 0$
and interpolates linearly between the values of $(x_n(h))$ at the
times $\{0,h,2h,\ldots\}$. As $h\to 0$, $\bar{x}_h$ approaches $x$
on any compact interval $[0,T]$ in the sense that
$\lim_{h\to0}\sup_{0\leq t\leq T} |x(t)-\bar{x}_h(t)|=0$ (see e.g.,
\cite{bellzenn:2003}).

\subsection{General discrete comparison results}
In this section we simply state our most general comparison results for the discretised equation. Later, we will apply these results to obtain concrete estimates of the growth of solutions of the discretised equation. 
\begin{theorem} \label{thm:limsup2disc}
Suppose that $f$ obeys (\ref{eq:f1}) and (\ref{eq:f2}).
Let $g$ be non--decreasing and obey (\ref{eq:g}) and let $\tau>0$
and $\psi\in C([-\tau,0];(0,\infty))$. Suppose that there exists a continuous function $\phi$ such that $\Gamma$, $\Gamma_c$ are defined by (\ref{def.Gamma}) and (\ref{def.Gammac}) respectively, and that $\Gamma$ obeys (\ref{eq.Gammatoinfty}). Suppose also that \eqref{eq.etaepstoeta}
%\begin{equation} \label{eq.etaepstoeta}
%\lim_{\epsilon\to 0^+}\eta(\varepsilon)=\eta,
%\end{equation}
and suppose that $f$ obeys \eqref{eq.fsmallg2},
and that $g$ and $\phi$ obey \eqref{eq.gphigamma} where $\bar{\eta}_\varepsilon$ obeys \eqref{eq.naretaeta}.
Suppose finally that $\phi$ and $f$ are non--decreasing. If $x_h$ is the unique solution of \eqref{eq:disceq}, then it obeys
\begin{equation}\label{eq.growthlimsupdisc}
\limsup_{n\to\infty} \frac{\Gamma(x_h(n))}{nh}\leq \eta.
\end{equation}
\end{theorem}

\begin{theorem} \label{thm:liminfdisc}
Suppose that $f$ obeys (\ref{eq:f1}) and (\ref{eq:f2}).
Let $g$ be non--decreasing and obey (\ref{eq:g}) and let $\tau>0$
and $\psi\in C([-\tau,0];(0,\infty))$. Suppose that there exists a continuous function $\phi$ such that $\Gamma$, $\Gamma_c$ are defined by (\ref{def.Gamma}) and (\ref{def.Gammac}) respectively, and $\Gamma$ obeys (\ref{eq.Gammatoinfty}). Suppose also that \eqref{eq.muepstomu} holds 
and that $g$ and $\phi$ obey
\begin{equation}  \label{eq.gphigamma2h}
\liminf_{x\to\infty} \frac{g(x)}{\phi(\Gamma_{\mu(\varepsilon)}^{-1}(\Gamma_{\mu(\varepsilon)}(x)+(\tau+h)(1-\epsilon)))}
=\bar{\mu}_\varepsilon\in (0,\infty] \quad \mbox{for every $\varepsilon\in (0,1)$},
\end{equation}
where \eqref{eq.narmumu} also holds. If $x_h$ is the unique solution of \eqref{eq:disceq}, then
\begin{equation} \label{eq.growthliminfdisc}
\liminf_{n\to\infty} \frac{\Gamma(x_h(n))}{nh}\geq \mu.
\end{equation}
\end{theorem}

\subsection{Preservation of growth rate for regularly varying $g$}
In \cite{AppMcCartRod09}, it was shown that the uniform Euler scheme \eqref{eq:disceq} and the continuous time extension $x_h$ preserves the rate of growth of the underlying continuous  equation \eqref{eq:cnseq} in the case when $g$ is in $\text{RV}_\infty(\beta)$ for $\beta\leq 1$, and $g$ is sublinear. We extract here the relevant parts of Theorems 2.4 and 2.5 of \cite{AppMcCartRod09}.
\begin{theorem}  \label{thm:presrv0}
Let $f$ obey (\ref{eq:f1}), (\ref{eq:f2}). Let $g$ obey (\ref{eq:g}). Let $\tau>0$ and $\psi\in C([-\tau,0];(0,\infty))$. Let $\beta\leq 1$ and suppose $g\in \mbox{RV}_\infty(\beta)$, and $\lim_{x\to 0}g(x)/x=0$. If $\lim_{x\to\infty} f(x)/g(x)=0$,
then the unique  solution $x_h$ of (\ref{eq:disceq}) obeys
\begin{equation} \label{eq.growthlimdelaydiscsublin}
%\frac{1}{\tau+h}\leq \liminf_{n\to\infty} \frac{G(x_h(n))}{nh}\leq
\lim_{n\to\infty} \frac{1}{nh}\int_1^{x_h(n)} \frac{1}{g(s)}\,ds=  1.
\end{equation}
Moreover, if $\bar{x}_h$ is the linear interpolant given by  \eqref{def.lininterpolant}, then
\begin{equation*} %\label{eq.growthlimdelaydiscsublininterp}
\lim_{t\to\infty} \frac{1}{t}\int_1^{\bar{x}_h(t)} \frac{1}{g(s)}\,ds=  1.
\end{equation*}
\end{theorem}

In this paper, we demonstrate that the essential growth rate is preserved for all $h>0$, and that the exact rate of
growth is recovered in the limit as $h\to 0^+$, in a sense now made precise. We first consider the discrete analogue of
Theorem~\ref{thm:limnoinst}.
\begin{theorem}  \label{thm:presrv1}
Let $f$ obey (\ref{eq:f1}), (\ref{eq:f2}). Let $g$ obey (\ref{eq:g}). Let $\tau>0$ and $\psi\in C([-\tau,0];(0,\infty))$. Suppose $g\in \mbox{RV}_\infty(1)$, $x\mapsto g(x)/x$ is asymptotic to a non--decreasing function, $\lim_{x\to\infty}g(x)/x=\infty$, and  $\lim_{x\to\infty}f(x)/x=0$.
If $G$ is defined by (\ref{def:G}), then the unique  solution $x_h$ of (\ref{eq:disceq}) obeys
\begin{equation} \label{eq.growthlimdelaydisc}
\lim_{n\to\infty} \frac{G(x_h(n))}{nh}=  \frac{1}{\tau+h}.
\end{equation}
Moreover, if $\bar{x}_h$ is the linear interpolant given by  \eqref{def.lininterpolant}, then
\begin{equation}  \label{eq.growthlimdelaydiscrv1interp}
\lim_{t\to\infty} \frac{G(\bar{x}_h(t))}{t}=\frac{1}{\tau+h}.
\end{equation}
\end{theorem}
The proof is postponed to the final section. By comparing \eqref{eq.growthlimrv1} and \eqref{eq.growthlimdelaydiscrv1interp}, it can be seen that the essential growth rate is recovered by the linear interpolant for all $h>0$, and the exact rate is recovered in the limit as $h\to 0^+$.

The rate of growth is also recovered in the same manner in the case when $g$ grows polynomially at a superlinear rate,
as confirmed by the following discrete analogue of Theorem~\ref{thm:polyg}.
\begin{theorem}  \label{thm:presrvbeta}
Let $f$ obey (\ref{eq:f1}), (\ref{eq:f2}). Let $g$ obey (\ref{eq:g}). Let $\tau>0$ and $\psi\in C([-\tau,0];(0,\infty))$. Suppose that there exists $\beta>1$ such that $g$ obeys
\[
\lim_{x\to\infty} \frac{\log g(x)}{\log x}=\beta,
\]
and $\lim_{x\to\infty}f(x)/x=0$.
Then the unique  solution $x_h$ of (\ref{eq:disceq}) obeys
\begin{equation} \label{eq.growthlimdelaydiscrvbeta}
\lim_{n\to\infty} \frac{\log_2 x_h(n)}{nh}= \frac{\log \beta}{\tau+h}.
\end{equation}
Moreover, if $\bar{x}_h$ is the linear interpolant given by  \eqref{def.lininterpolant}, then
\begin{equation}   \label{eq.growthlimdelaydiscrvbetainterp}
\lim_{t\to\infty}\frac{\log_2 \bar{x}_h(t)}{t}=\frac{\log \beta}{\tau+h}.
\end{equation}
\end{theorem}
Once again, by comparing \eqref{eq.polygrowth} and \eqref{eq.growthlimdelaydiscrvbetainterp}, we see that the essential growth rate is recovered by the linear interpolant for all $h>0$, and the exact rate is recovered in the limit as $h\to 0^+$. Again, we relegate the proof to the end.

\section{Proof of Main Continuous--Time Results} \label{sec.continuousproofs}
In this section, we give the proofs of the main results from Section~\ref{sec.continuousgrowth}, with the exception of Theorem~\ref{thm:limnoinst}, whose proof is strongly based on that of Theorem~\ref{thm:presrv1}. The proofs of these two results, along with Theorem~\ref{thm:presrvbeta}, are given in Section~\ref{sec.discreteproofs}.
\subsection{Proof of Theorem~\ref{thm:thmcnsnonexplosion}}
Suppose that $x$ has a finite interval of existence. Then there is a unique continuous solution of
\eqref{eq:cnseq} on $[-\tau,T)$ where $T\in (0,\infty]$ is such
that
\[
\lim_{t\to T^-} x(t)=\infty.
\]
The limit is $+\infty$ because the positivity of the initial condition, together with the non--negativity of
$f$ and $g$ ensure that $x'(t)\geq 0$ for $t\in [0,T)$.

We wish to rule out the possibility that $T<+\infty$. Suppose that
$T\in(0,\tau]$. Clearly, if $g_1=\max_{s\in[-\tau,0]} g(x(s))\geq
0$, we have
\[
x'(t)\leq f(x(t))+g_1, \quad t\in [0,T).
\]
Define $f_1(x):=f(x)+g_1$ for $x\geq 0$. Then, as $x(t)\to\infty$
as $t\to T^-$, we have
\[
\int_{x(0)}^\infty \frac{1}{f_1(x)}\,dx=\lim_{t\to T^-}\int_0^t
\frac{x'(s)}{f_1(x(s))}\,ds \leq T.
\]
However, \eqref{eq:f2} implies that $\int_{x(0)}^\infty
1/f_1(u)\,du=\infty$, which gives a contradiction. Hence $T>\tau$.

Suppose now that $x$ does not explode in $[0,n\tau]$, but does in
$(n\tau,(n+1)\tau]$. This is true for $n=1$. Clearly, if
$g_n=\max_{s\in[(n-1)\tau,n\tau]} g(x(s))\geq 0$, we have
\[
x'(t)\leq f(x(t))+g_n, \quad t\in [n\tau,T).
\]
Define $f_n(x):=f(x)+g_n$ for $x\geq 0$. Then, as $x(t)\to\infty$
as $t\to T^-$, we have
\[
\int_{x(n\tau)}^\infty \frac{1}{f_n(x)}\,dx=\lim_{t\to
T^-}\int_{n\tau}^t \frac{x'(s)}{f_n(x(s))}\,ds \leq T-n\tau.
\]
However, \eqref{eq:f2} implies that $\int_{x(n\tau)}^\infty
1/f_n(u)\,du=\infty$, which gives a contradiction. Hence
$T>(n+1)\tau$. Since this is true for any $n\in\mathbb{N}$, it
follows that $T=\infty$.

We have shown that \eqref{eq:cnseq} has interval of existence $[-\tau,\infty)$. Since $\psi(t)>0$ for $t\in [-\tau,0]$
and $f(x)\geq 0$, $g(x)\geq 0$ for all $x\geq 0$, we have that $x'(t)\geq 0$ for all $t\geq 0$. Therefore
$\lim_{t\to\infty} x(t)=:L\in [\psi(0),\infty]$. Suppose that $L>0$ is finite. Since
\[
x(t)=\psi(0)+\int_0^t f(x(s))\,ds + \int_0^{t-\tau} g(x(s))\,ds, \quad t\geq \tau,
\]
by the continuity of $f$ and $g$ we have
\[
\lim_{t\to\infty}\frac{1}{t}\int_0^t f(x(s))\,ds=f(L), \quad
\lim_{t\to\infty}\frac{1}{t}\int_0^{t-\tau} g(x(s))\,ds = g(L).
\]
Since $x(t)$ tends to the finite limit $L$, we get
\[
0=\lim_{t\to\infty}
\frac{x(t)}{t}=
\lim_{t\to\infty}\frac{\psi(0)}{t}+\frac{1}{t}\int_0^t f(x(s))\,ds +
\frac{1}{t}\int_0^{t-\tau} g(x(s))\,ds = f(L)+g(L).
\]
Since $g$ is positive and $f$ is nonnegative, we have $L=0$, a contradiction. Hence $x$ obeys (\ref{eq:cnsnonexplode}),
as claimed.

\subsection{Proof of Theorem~\ref{thm:limsup2}}
By (\ref{eq.gphigamma}) for every $\epsilon\in (0,1)$ there exists $x_2(\epsilon)>0$ such that for $x>x_2(\epsilon)$ we have
\[
g(x)<(\bar{\eta}_\epsilon+\epsilon)\phi(\Gamma_{\eta(\epsilon)}^{-1}(\Gamma_{\eta(\epsilon)}(x)+\tau))
\leq (\bar{\eta}+\epsilon)\phi(\Gamma_{\eta(\epsilon)}^{-1}(\Gamma_{\eta(\epsilon)}(x)+\tau)),
\]
where the last inequality is a consequence of (\ref{eq.naretaeta}). Since $\bar{\eta}<\eta=\lim_{\epsilon\rightarrow 0^+} \eta(\epsilon)$,
there exists $\epsilon'\in (0,1)$ such that for $\epsilon<\epsilon'$, we have $\eta(\epsilon)>\bar{\eta}+\epsilon$. Thus for all $\epsilon<\epsilon'<1$ we have
\begin{equation} \label{eq.gltphi}
g(x)<\eta(\epsilon)\phi(\Gamma_{\eta(\epsilon)}^{-1}(\Gamma_{\eta(\epsilon)}(x)+\tau)), \quad x>x_2(\epsilon).
\end{equation}
By (\ref{eq.fsmallg2}) for every $\epsilon\in (0,1)$ there exists an $x_1(\epsilon)>0$ such that
\begin{equation} \label{eq.fltg}
f(x)\leq \epsilon \eta(\epsilon) \phi(x), \quad x>x_1(\epsilon).
\end{equation}
%Define $y_2(\epsilon)=\Gamma_{\eta(\epsilon)}(x_2(\epsilon))+\tau$. Now, if $x>x_2(\epsilon)$ we have that $y:=\Gamma_{\eta(\epsilon)}(x)+\tau$
%obeys $y>y_2(\epsilon)$ and $x=\Gamma_{\eta(\epsilon)}^{-1}(y-\tau)$. Thus for $y>y_2(\epsilon)$, by (\ref{eq.gltphi}) we have
%$g(\Gamma_{\eta(\epsilon)}^{-1}(y-\tau))<\eta(\epsilon)\phi(\Gamma_{\eta(\epsilon)}^{-1}(y))$ for $y>y_2(\epsilon)$. Therefore
%\begin{equation} \label{eq.gltphi2}
%(1+\epsilon)g(\Gamma_{\eta(\epsilon)}^{-1}(x-\tau))<\eta(\epsilon)(1+\epsilon)\phi(\Gamma_{\eta(\epsilon)}^{-1}(x)), \quad x>y_2(\epsilon).
%\end{equation}
Define
\begin{equation} \label{def.c}
c(\epsilon)=\Gamma_{\eta(\epsilon)}(\psi^\ast+x_1(\epsilon)+x_2(\epsilon))+(1+\epsilon)\tau,
\end{equation}
and define also
\begin{equation}\label{def.xepsupper}
x_\epsilon(t)=\Gamma_{\eta(\epsilon)}^{-1}((1+\epsilon)t+c(\epsilon)),\quad  t\geq -\tau.
\end{equation}
This function is well--defined since $c(\epsilon)>\Gamma_{\eta(\epsilon)}(\psi^\ast)+(1+\epsilon)\tau$, so $c(\epsilon)-(1+\epsilon)\tau>\Gamma_{\eta(\epsilon)}(\psi^\ast)$, or $x_{\epsilon}(t)>\psi^\ast$ for all $t\in [-\tau,0]$.
Since $c(\epsilon)>\Gamma_{\eta(\epsilon)}(x_1(\epsilon))+(1+\epsilon)\tau$ and $\Gamma_{\eta(\epsilon)}$ is increasing,
$\Gamma_{\eta(\epsilon)}^{-1}(c(\epsilon)-(1+\epsilon)\tau)>x_1(\epsilon)$, so $x_\epsilon(t)>x_1(\epsilon)$ for all $t\geq -\tau$. Therefore by (\ref{eq.fltg}), $f(x_\epsilon(t))\leq \epsilon\eta(\epsilon) \phi(x_\epsilon(t))$. %g(\Gamma_{\eta(\varepsilon)}^{-1}(\Gamma_{\eta(\varepsilon)}(x_\epsilon(t))-\tau))
%=\epsilon g(\Gamma_{\eta(\epsilon)}^{-1}((1+\epsilon)t+c(\epsilon)-\tau))$. 
Also for $t\geq 0$, we have
\begin{align*}
g(x_\epsilon(t-\tau))
&=g(\Gamma_{\eta(\epsilon)}^{-1}((1+\epsilon)(t-\tau)+c(\epsilon))
=g(\Gamma_{\eta(\epsilon)}^{-1}((1+\epsilon)t-\tau-\epsilon \tau + c(\epsilon)))\\
&<g(\Gamma_{\eta(\epsilon)}^{-1}((1+\epsilon)t-\tau+c(\epsilon))).
\end{align*}
Now, because $c(\epsilon)>\Gamma_{\eta(\epsilon)}(x_2(\epsilon))+\tau$, we have that the argument of $g$ on 
the righthand side exceeds $x_2(\epsilon)$ for all $t\geq 0$. Therefore by \eqref{eq.gltphi}, we have 
\begin{align*}
g(x_\epsilon(t-\tau))&<g(\Gamma_{\eta(\epsilon)}^{-1}((1+\epsilon)t-\tau+c(\epsilon)))\\
&<\eta(\epsilon)\phi(
\Gamma_{\eta(\epsilon)}^{-1}(\Gamma_{\eta(\epsilon)}(\Gamma_{\eta(\epsilon)}^{-1}((1+\epsilon)t-\tau+c(\epsilon)))+\tau)
)\\
&=\eta(\epsilon)\phi(\Gamma_{\eta(\epsilon)}^{-1}((1+\epsilon)t-\tau+c(\epsilon))+\tau)\\
&=\eta(\epsilon)\phi(\Gamma_{\eta(\epsilon)}^{-1}((1+\epsilon)t+c(\epsilon)))\\
&=\eta(\epsilon)\phi(x_\epsilon(t)).
\end{align*}
Hence for $t\geq 0$
\begin{equation} \label{eq.fplusg}
f(x_\epsilon(t))+g(x_\epsilon(t-\tau))<(1+\epsilon)\eta(\epsilon) \phi(x_\epsilon(t)).
\end{equation}
Now for $t>0$, $\Gamma_{\eta(\epsilon)}(x_\epsilon(t))=(1+\epsilon)t + c(\epsilon)$, so
$\Gamma_{\eta(\epsilon)}'(x_\epsilon(t)) x_\epsilon'(t)=(1+\epsilon)$, or
$x_\epsilon'(t)=(1+\epsilon)\eta(\epsilon)\phi(x_\epsilon(t))$. Hence
\begin{equation} \label{eq.xepsdiffeq1}
x_\epsilon'(t)=(1+\epsilon)
\eta(\epsilon)\phi(\Gamma_{\eta(\epsilon)}^{-1}((1+\epsilon)t+c(\epsilon))), \quad t>0.
\end{equation}
%Now $c(\epsilon)>\Gamma_{\eta(\epsilon)}(x_2(\epsilon))+\tau=y_2(\epsilon)$. Then, because $(1+\epsilon)t+c(\epsilon)\geq c(\epsilon)$ for $t\geq 0$, we have
%$(1+\epsilon)t+c(\epsilon)>y_2(\epsilon)$ for $t\geq 0$. Thus by
%(\ref{eq.gltphi2}) we have
%$(1+\epsilon)g(\Gamma_{\eta(\epsilon)}^{-1}((1+\epsilon)t+c(\epsilon)-\tau))
%<\eta(\epsilon)(1+\epsilon)\phi(\Gamma_{\eta(\epsilon)}^{-1}((1+\epsilon)t+c(\epsilon)))$,
Thus by (\ref{eq.fplusg}) and (\ref{eq.xepsdiffeq1}) for $t>0$ we have
$x_\epsilon'(t)>f(x_\epsilon(t))+g(x_\epsilon(t-\tau))$.
%\eta(\epsilon)(1+\epsilon)\phi(\Gamma_{\eta(\epsilon)}^{-1}((1+\epsilon)t+c(\epsilon))) =x_\epsilon'(t)$.

Now as $x_\epsilon(t)>\psi^\ast=\max_{t\in [-\tau,0]}\psi(s)$, we have $x_\epsilon(t)>x(t)$ for $t\in [-\tau,0]$ and $x_\epsilon'(t)>f(x_\epsilon(t))+g(x_\epsilon(t-\tau))$ for $t\geq 0$.
Suppose that there is a $t_0>0$ such that $x_\epsilon(t)>x(t)$ for $t\in [-\tau,t_0)$
$x_\epsilon(t_0)=x(t_0)$. Therefore $x_\epsilon'(t_0)\leq x'(t_0)$. Then as
$g$ is non--decreasing,
\begin{align*}
x_\epsilon'(t_0)&\leq x'(t_0)=f(x(t_0))+g(x(t_0-\tau))\\
&=f(x_\epsilon(t_0))+g(x(t_0-\tau))\leq f(x_\epsilon(t_0))+g(x_\epsilon(t_0-\tau))\\
&<x_\epsilon'(t_0),
\end{align*}
a contradiction. Thus $x_\epsilon(t)>x(t)$ for all $t\geq -\tau$. Hence
$\Gamma_{\eta(\epsilon)}(x(t))<\Gamma_{\eta(\epsilon)}(x_\epsilon(t))$
for all $t\geq -\tau$. Hence
\[
\Gamma_{\eta(\epsilon)}(x(t))<\Gamma_{\eta(\epsilon)}(x_\epsilon(t))=(1+\epsilon)t + c(\epsilon), \quad t \geq -\tau.
\]
But $\Gamma(x(t))=\eta(\epsilon)\Gamma_{\eta(\epsilon)}(x(t))
<(1+\epsilon)\eta(\epsilon)t+\eta(\epsilon)c(\epsilon)$. Therefore
\[
\limsup_{t\rightarrow\infty} \Gamma(x(t))/t\leq (1+\epsilon)\eta(\epsilon).
\]
Since $\epsilon>0$ is arbitrary, and $\eta(\epsilon)\rightarrow\eta$ as $\epsilon\rightarrow 0$, we have (\ref{eq.growthlimsup}).

\subsection{Proof of Theorem~\ref{thm:rvbetalt1}}
Suppose that $\phi(x)=g(x)$ for $x>0$. Thus $\Gamma_\eta(x)=\eta^{-1}\int_{\psi^\ast}^x du/g(u)$. Let $z(t)=\Gamma_\eta^{-1}(t)$ for $t\geq 0$. Then $z'(t)=\eta g(z(t))$ for $t>0$ with $z(0)=\psi^\ast$. Thus $z'(t)/z(t)\to 0$ as $t\to\infty$. Therefore
\[
\log\left(\frac{z(t)}{z(t-\theta)}\right)=\int_{t-\theta}^t \frac{z'(s)}{z(s)}\,ds \to 0 \quad\mbox{as $t\to\infty$},
\]
so $\lim_{t\to\infty} z(t-\theta)/z(t)=1$ for any $\theta\in\mathbb{R}$. Since $g\in\mbox{RV}_\infty(\beta)$, we have
\[
\lim_{t\to\infty} g(z(t-\theta))/g(z(t))=1.
\]
Hence $\lim_{t\to\infty} g(\Gamma_\eta^{-1}(t-\theta))/g(\Gamma_\eta^{-1}(t))=1$. Since $\Gamma_\eta^{-1}(t)\to \infty$ as $t\to\infty$, we have
\begin{equation} \label{eq.ggarglim1}
\lim_{x\to\infty} \frac{g(x)}{\phi(\Gamma_\eta^{-1}(\Gamma_\eta(x)+\theta))}
=
\lim_{x\to\infty} \frac{g(x)}{g(\Gamma_\eta^{-1}(\Gamma_\eta(x)+\theta))}=1.
\end{equation}
Since this holds for every $\eta>0$ and $\theta\in\mathbb{R}$ it follows that (\ref{eq.gphigamma}) and (\ref{eq.gphigamma2}) hold with $\bar{\eta}_\epsilon=\bar{\mu}_\epsilon=1$.
Let $\rho\in (0,1)$. Define $\mu(\epsilon)=1-\rho$ and $\eta(\epsilon)=1+\rho$. Then with $\eta=1+\rho$
and $\mu=1-\rho$, (\ref{eq.etaepstoeta}), (\ref{eq.muepstomu}), (\ref{eq.naretaeta}) and (\ref{eq.narmumu}) hold.
To prove (\ref{eq.fsmallg2}), we note that
\[
\lim_{x\to\infty} \frac{f(x)}{\phi(x)}=\lim_{x\to\infty} \frac{f(x)}{g(x)}=0.
\]
%
%\[
%\lim_{x\to\infty} \frac{f(x)}{g(\Gamma_{\eta}^{-1}(\Gamma_{\eta}(x)-\tau))}
%=
%\lim_{x\to\infty} \frac{f(x)}{g(x)}\cdot \frac{g(x)}{g(\Gamma_{\eta}^{-1}(\Gamma_{\eta}(x)-\tau))}=0,
%\]
%because (\ref{eq.ggarglim1}) holds and $f(x)/g(x)\to 0$ as $x\to\infty$.
Since all the hypotheses of Theorems~\ref{thm:limsup} and \ref{thm:liminf} hold, we have $\limsup_{t\to\infty} \Gamma(x(t))/t\leq 1+\rho$ and $\liminf_{t\to\infty} \Gamma(x(t))/t\geq 1-\rho$. Letting $\rho\to 0$, we have $\lim_{t\to\infty} \Gamma(x(t))/t=1$, whence the result.

\subsection{Proof of  Theorem~\ref{thm:limrv1}}
Since $g\in \mbox{RV}_\infty(1)$, it follows that there exists an increasing and continuously differentiable function $\delta:[\psi^\ast,\infty)\to (0,\infty)$ with $\delta(\psi^\ast)>e\psi^\ast$ such that $\delta(x)/g(x)\to 1$ as $x\to\infty$ and $x\delta'(x)/\delta(x)\to 1$ as $x\to\infty$. Define $\phi(x)=x\log(\delta(x)/x)$ for $x\geq \psi^\ast$. Define $\Gamma(x)=\int_{\psi^\ast}^x du/\phi(u)$ for $x\geq \psi^\ast$. Since $(g(x)/x)/(\delta(x)/x)\to 1$ as $x\to\infty$, we have $\log(g(x)/x)/\log(\delta(x)/x)\to 1$ as $x\to\infty$. Therefore by L'H\^opital's rule, we have $\Gamma(x)/G(x)\to 1$ as $x\to\infty$.

Define $\Gamma_\eta(x)=\Gamma(x)/\eta$ and $\delta_1(x)=\delta(x)/x$ for $x\geq \psi^\ast$. Since
$x\delta'(x)/\delta(x)\to 1$ as $x\to\infty$, we have that $\delta_1$ is continuously differentiable and  $x\delta_1'(x)/\delta_1(x)\to 0$ as $x\to\infty$. Define $y(t)=\Gamma_{\eta}^{-1}(t)$ for $t\geq 0$ and $u(t)=\log \delta_1(y(t))$. Then $y'(t)=\eta \phi(y(t))=\eta y(t)\log\delta_1(y(t))=\eta y(t) u(t)$. Moreover
since $\Gamma_\eta(x)\to\infty$ as $x\to\infty$, we have that $y(t)\to\infty$ as $t\to\infty$.
Thus
\begin{align*}
\lim_{x\to\infty} \frac{g(x)}{\phi(\Gamma_\eta^{-1}(\Gamma_\eta(x)+\theta))}
&=\lim_{x\to\infty} \frac{\delta(x)}{\phi(\Gamma_\eta^{-1}(\Gamma_\eta(x)+\theta))}\\
&=\lim_{t\to\infty} \frac{\delta(\Gamma_\eta^{-1}(t-\theta))}{\Gamma_{\eta}^{-1}(t)\log(\delta(\Gamma_\eta^{-1}(t))/t)},
\end{align*}
and therefore we have
\[
\lim_{x\to\infty} \frac{g(x)}{\phi(\Gamma_\eta^{-1}(\Gamma_\eta(x)+\theta))}
=\lim_{t\to\infty} \frac{y(t-\theta)\delta_1(y(t-\theta))}{y(t)\log \delta_1(y(t))}.
\]
Since $\log(y(t)/y(t-\theta))=\int_{t-\theta}^t y'(s)/y(s)\,ds = \int_{t-\theta}^t \eta u(s)\,ds$. Hence
\begin{eqnarray*}
\lefteqn{\log\left(
\lim_{x\to\infty} \frac{g(x)}{\phi(\Gamma_\eta^{-1}(\Gamma_\eta(x)+\theta))}\right)}
\\
&=&\lim_{t\to\infty} \{-\log (y(t)/y(t-\theta)) + u(t-\theta)  - \log u(t)\}\\
&=&\lim_{t\to\infty} u(t)\left\{-\eta\frac{1}{u(t)}\int_{t-\theta}^t u(s)\,ds + \frac{u(t-\theta)}{u(t)}  - \frac{\log u(t)}{u(t)}\right\}.
\end{eqnarray*}
Since $\delta_1$, $y$ are continuously differentiable, so is $u$, and we have
\[
u'(t)=\delta_1'(y(t))y'(t)/\delta_1(y(t))=\eta u(t) \cdot y(t)\delta_1'(y(t))/\delta_1(y(t)).
\]
Since $x\delta_1'(x)/\delta_1(x)\to 0$ as $x\to\infty$ and $y(t)\to\infty$ as $t\to\infty$, we have $u'(t)/u(t)\to 0$ as $t\to\infty$. Also we have $u(t)\to\infty$ as $t\to\infty$. Therefore $u(t-\theta)/u(t)\to 1$ as $t\to\infty$ and
\[
\lim_{t\to\infty}\int_{t-\theta}^t u(s)\,ds/u(t)= \theta,
\]
so
\[
\lim_{t\to\infty}
\left\{-\eta\frac{1}{u(t)}\int_{t-\theta}^t u(s)\,ds + \frac{u(t-\theta)}{u(t)}  - \frac{\log u(t)}{u(t)}\right\}
=1-\eta \theta.
\]
Therefore we have
\[
\log\left(
\lim_{x\to\infty} \frac{g(x)}{\phi(\Gamma_\eta^{-1}(\Gamma_\eta(x)+\theta))}\right)
=\left\{
\begin{array}{cc}
-\infty & \mbox{if $1-\eta\theta<0$} \\
+\infty & \mbox{if $1-\eta\theta>0$}.
\end{array}
\right.
%\lim_{t\to\infty} u(t)\left\{-\eta\frac{1}{u(t)}\int_{t-\theta}^t u(s)\,ds + \frac{u(t-\theta)}{u(t)}  - \frac{\log u(t)}{u(t)}\right\}.
\]
Therefore, with $\eta(\epsilon)=(1+\epsilon)/\tau$ and $\mu(\epsilon)=(1-\epsilon)/\tau$, we have
\[
\lim_{x\to\infty} \frac{g(x)}{\phi(\Gamma_{\eta(\epsilon)}^{-1}(\Gamma_{\eta(\epsilon)}(x)+\tau))}=0, \quad
\lim_{x\to\infty} \frac{g(x)}{\phi(\Gamma_{\mu(\epsilon)}^{-1}(\Gamma_{\mu(\epsilon)}(x)+\tau(1-\epsilon)))}=\infty.
\]
Since $\mu(\epsilon),\eta(\epsilon)\to 1/\tau$ as $\epsilon\to 0$, and we have $\bar{\eta}_\epsilon=0=:\bar{\eta}<1/\tau$
and $\bar{\mu}_\epsilon=+\infty=:\bar{\mu}>1/\tau$.

Next, note that \eqref{eq.fsmallvsgrv1} implies 
\[
\lim_{x\to\infty} \frac{f(x)}{\phi(x)}
=\lim_{x\to\infty} \frac{f(x)}{x\log(\delta(x)/x)}
=\lim_{x\to\infty} \frac{f(x)}{x\log(g(x)/x)}=0.
\]
%so $\lim_{x\to\infty} f(x)/g(\Gamma_{\eta}^{-1}(\Gamma_\eta(x)-\tau))=
%\lim_{t\to\infty} f(y(t+\tau))/\delta(y(t))$,
%so
%\begin{multline*}
%\log\left( \lim_{x\to\infty}
% \frac{f(x)}{g(\Gamma_{\eta}^{-1}(\Gamma_\eta(x)-\tau))}\right)
% =\lim_{t\to\infty} \Biggl\{\left(\frac{\log f(y(t+\tau))}{\log y(t+\tau)}-1\right)\log y(t+\tau)
%\\ +\log y(t+\tau)-\log y(t)-u(t)\Biggr\}.
%\end{multline*}
%Since $\log y(t+\tau)-\log y(t)=\int_t^{t+\tau} \eta u(s)\,ds$, we have
%\begin{eqnarray*}
%\lefteqn{\log\left( \lim_{x\to\infty}
% \frac{f(x)}{g(\Gamma_{\eta}^{-1}(\Gamma_\eta(x)-\tau))}\right)}\\
% &=&\lim_{t\to\infty} \log y(t+\tau)\Biggl\{\left(\frac{\log f(y(t+\tau))}{\log y(t+\tau)}-1\right)
% \\
% &&\qquad+\frac{u(t+\tau)}{\log y(t+\tau)}
% \left(\frac{1}{u(t+\tau)}\eta\int_t^{t+\tau} u(s)\,ds-\frac{u(t)}{u(t+\tau)}\right)\Biggr\}.
%\end{eqnarray*}
%Properties of $u$ imply that $\eta\int_t^{t+\tau} u(s)\,ds/u(t+\tau)-u(t)/u(t+\tau)\to \eta\tau-1$ as $t\to\infty$.
%Since $\delta_1\in \mbox{RV}_\infty(0)$, we have $\log \delta_1(x)/\log x\to 0$ as $x\to\infty$.
%Therefore $u(t)/\log y(t)=\log \delta_1(y(t))/\log y(t)\to 0$ as $t\to\infty$. Hence the quantity inside the
%curly brackets has limit superior not exceeding $\gamma-1$, which is negative. Since $y(t+\tau)\to \infty$ as $t\to\infty$ we have that $\lim_{x\to\infty} f(x)/g(\Gamma_{\eta(\epsilon)}^{-1}(\Gamma_{\eta(\epsilon)}(x)-\tau))=0$.
Therefore by Theorem~\ref{thm:lim}, we have $\lim_{t\to\infty} \Gamma(x(t))/t=1/\tau$, and due to the fact that  $\lim_{x\to\infty}G(x)/\Gamma(x)=1$, we get $\lim_{t\to\infty} G(x(t))/t=1/\tau$, as required.

\subsection{Proof of Theorem~\ref{ex:4}}
Set $\phi(x)=x$ for $x\geq \psi^\ast$. Then $\Gamma_\eta(x)=\eta^{-1}\log(x/\psi^\ast)$, $\Gamma_\eta^{-1}(x)=\psi^\ast e^{\eta x}$,
and $\phi(\Gamma_\eta^{-1}(\Gamma_\eta(x)+\theta))=xe^{\eta \theta}$. Thus $\lim_{x\to\infty} g(x)/\phi(\Gamma_\eta^{-1}(\Gamma_\eta(x)+\theta))=Ce^{-\eta\theta}$.
Define $c(\nu):=\nu-Ce^{-\nu \tau}$. Then $c$ is increasing on $[0,\infty)$ and there is a unique $\lambda>0$
such that $c(\lambda)=0$, or $\lambda=Ce^{-\lambda \tau}$. Let $\sigma\in\mathbb{R}$ and $\lambda_\sigma:=\lambda(1+\sigma)$. For $\sigma>0$, $c(\lambda_\sigma)>0$ or
$\lambda_\sigma>Ce^{-\lambda_\sigma \tau}$. Similarly, $\lambda_{-\sigma}<Ce^{-\lambda_{-\sigma} \tau}$. Define $\eta(\epsilon)=\lambda_\sigma(1+\epsilon)$. Then $\eta(\epsilon)\to \lambda_\sigma=:\eta$ as $\epsilon\to 0$. Also
$\lim_{x\to\infty} g(x)/\phi(\Gamma_{\eta(\epsilon)}^{-1}(\Gamma_{\eta(\epsilon)}(x)+\tau))
=Ce^{-\lambda_\sigma(1+\epsilon)\tau}=:\bar{\eta}_\epsilon$. Then $\sup_{\epsilon\in (0,1)} \bar{\eta}_\epsilon=C e^{-\lambda_\sigma\tau}=:\bar{\eta}$. But $\bar{\eta}=C e^{-\lambda_\sigma\tau}<\lambda_\sigma=\eta$.
Finally, $f(x)/\phi(x)=f(x)/x\to 0$ as $x\to\infty$, 
%$g(\Gamma_{\eta(\epsilon)}^{-1}(\Gamma_{\eta(\epsilon)}(x)-\tau))/x\to Ce^{-\eta(\epsilon)\tau}$ as $x\to\infty$. Thus 
%$f(x)/g(\Gamma_{\eta(\epsilon)}^{-1}(\Gamma_{\eta(\epsilon)}(x)-\tau))\to 0$ as $x\to\infty$, 
and so by Theorem~\ref{thm:limsup}, $\limsup_{t\to\infty} \Gamma(x(t))/t\leq \lambda_\sigma$, or $\limsup_{t\to\infty} \log x(t)/t\leq \lambda(1+\sigma)$. Letting $\sigma\downarrow 0$ yields $\limsup_{t\to\infty} \log x(t)/t\leq \lambda$. Define $\mu(\epsilon)=\lambda_{-\sigma}(1-\epsilon)$. Then $\lim_{\epsilon\to 0}\mu(\epsilon)=\lambda_{-\sigma}=:\mu$.
Also $\lim_{x\to\infty} g(x)/\phi(\Gamma_{\mu(\epsilon)}^{-1}(\Gamma_{\mu(\epsilon)}(x)+\tau(1-\epsilon)))
=Ce^{-\lambda_{-\sigma}(1-\epsilon)\tau}=:\bar{\mu}_\epsilon$.
Then $\inf_{\epsilon\in (0,1)} \bar{\mu}_\epsilon=C e^{-\lambda_{-\sigma}\tau}=:\bar{\mu}$. But
$\bar{\mu}=C e^{-\lambda_{-\sigma}\tau}>\lambda_{-\sigma}=\mu$.
Thus by Theorem~\ref{thm:liminf}, $\liminf_{t\to\infty} \Gamma(x(t))/t\geq \lambda_{-\sigma}$, or $\liminf_{t\to\infty} \log x(t)/t\geq \lambda(1-\sigma)$. Letting $\sigma\downarrow 0$ yields $\liminf_{t\to\infty} \log x(t)/t\geq \lambda$, whence the result.
%Combining the results gives $\lim_{t\to\infty} \log x(t)/t=\lambda$.
%\sup_{\epsilon\in (0,1)} e^{-\lambda_\sigma\epsilon\tau}$

\subsection{Proof of Theorem~\ref{thm:polyg}}
Define $\phi(x)=(1+x)\log(1+x)$ for $x\geq \psi^\ast$. Hence for $\eta>0$ we have
\[
\Gamma_\eta(x)=\frac{1}{\eta}\log\left(\frac{\log (1+x)}{\log (1+\psi^\ast)}  \right),
\quad
\Gamma_\eta^{-1}(x)=\exp\left(\log(1+\psi^\ast)e^{\eta x}\right)-1.
\]
Thus $\phi(\Gamma_\eta^{-1}(\Gamma_\eta(x)+\theta))=e^{\eta\theta}(1+x)^{e^{\eta \theta}} \log(1+x)$.
Also $\Gamma_\eta^{-1}(\Gamma_\eta(x)-\tau)=(1+x)^{e^{-\eta \tau}}-1$.
Therefore
\begin{eqnarray*}
\lim_{x\to\infty} \frac{g(x)}{\phi(\Gamma_\eta^{-1}(\Gamma_\eta(x)+\theta))}
&=&e^{-\eta\theta}\lim_{x\to\infty} \frac{g(x)}{(1+x)^{e^{\eta \theta}} \log(1+x)}, \\
\lim_{x\to\infty} \frac{f(x)}{g(\Gamma_\eta^{-1}(\Gamma_\eta(x)-\tau))}
&=&\lim_{x\to\infty} \frac{f(x)}{g((1+x)^{e^{-\eta \tau}}-1)}.
\end{eqnarray*}
Next, $\eta(\epsilon):=\epsilon+\log(\beta)/\tau$. Then $\lim_{\epsilon\to 0}\eta(\epsilon)=\log(\beta)/\tau=:\eta$, and so
\[
\lim_{x\to\infty} g(x)/\phi(\Gamma_{\eta(\epsilon)}^{-1}(\Gamma_{\eta(\epsilon)}(x)+\tau))=0.
\]
Therefore $\bar{\eta}_\epsilon=0$, so $\bar{\eta}=0<\log(\beta)/\tau=\eta$.
Next, as $f(x)/(x\log x)\to 0$ as $x\to\infty$, we have
\[
\lim_{x\to\infty} \frac{f(x)}{\phi(x)}=\frac{f(x)}{(1+x)\log(1+x)}=0.
\]
%Also, as $\gamma- \beta e^{-\eta(\epsilon)\tau}=\gamma- e^{-\tau\epsilon}<0$ for all $\epsilon$ sufficiently small,
%\begin{align*}
%\lefteqn{\log\left(\lim_{x\to\infty} \frac{f(x)}{g(\Gamma_{\eta(\epsilon)}^{-1}(\Gamma_{\eta(\epsilon)}(x)-\tau))}\right)}\\
%&=\lim_{x\to\infty} \log x\left(\frac{\log f(x)}{\log x}
%- \frac{\log g((1+x)^{e^{-\eta(\epsilon)\tau}}-1)}{\log((1+x)^{e^{-\eta(\epsilon)\tau}}-1)}\cdot \frac{\log((1+x)^{e^{-\eta(\epsilon)\tau}}-1)}{\log x}\right)
%\\
%&=-\infty,
%\end{align*}
%and so $\lim_{x\to\infty} f(x)/g(\Gamma_{\eta(\epsilon)}^{-1}(\Gamma_{\eta(\epsilon)}(x)-\tau))=0$.
By Theorem~\ref{thm:limsup2},
\[
\limsup_{t\to\infty} \Gamma(x(t))/t\leq \eta,
\]
or equivalently $\limsup_{t\to\infty} \log\log x(t)/t\leq \log(\beta)/\tau$.
We now obtain a lower bound. Define $\mu(\epsilon)=\log(\beta)/\tau$ for $\epsilon>0$. Then
\[
\lim_{x\to\infty} \frac{g(x)}{\phi(\Gamma_{\mu(\epsilon)}^{-1}(\Gamma_{\mu(\epsilon)}(x)+\tau(1-\epsilon)))}
=\beta^{-(1-\epsilon)}\lim_{x\to\infty} \frac{g(x)}{(1+x)^{\beta^{1-\epsilon}} \log(1+x)}.
\]
Therefore
\[
\lim_{x\to\infty} g(x)/\phi(\Gamma_{\mu(\epsilon)}^{-1}(\Gamma_{\mu(\epsilon)}(x)+\tau(1-\epsilon)))=\infty
\],
so $\bar{\mu}_\epsilon=+\infty=\bar{\mu}>\mu=\log(\beta)/\tau$. By Theorem~\ref{thm:liminf},
$\liminf_{t\to\infty} \Gamma(x(t))/t\geq \mu$, or
$\liminf_{t\to\infty} \log\log x(t)/t\geq \log \beta/\tau$, which proves (\ref{eq.polygrowth}).

\section{Proof of Main Discrete--Time Results} \label{sec.discreteproofs}
In this section, we give the proofs of results from Section \ref{sec.discretegrowth}. We also give the proof of Theorem~\ref{thm:limnoinst}, which is greatly facilitated by the proof of Theorem~\ref{thm:presrv1}.
\subsection{Proof of Theorem~\ref{thm:limsup2disc}}
By (\ref{eq.gphigamma}) for every $\epsilon\in (0,1)$ there exists $x_2(\epsilon)>0$ such that for $x>x_2(\epsilon)$ we have
\[
g(x)<(\bar{\eta}_\epsilon+\epsilon)\phi(\Gamma_{\eta(\epsilon)}^{-1}(\Gamma_{\eta(\epsilon)}(x)+\tau))
\leq (\bar{\eta}+\epsilon)\phi(\Gamma_{\eta(\epsilon)}^{-1}(\Gamma_{\eta(\epsilon)}(x)+\tau)),
\]
where the last inequality is a consequence of (\ref{eq.naretaeta}). Since $\bar{\eta}<\eta=\lim_{\epsilon\rightarrow 0^+} \eta(\epsilon)$,
there exists $\epsilon'\in (0,1)$ such that for $\epsilon<\epsilon'$, we have $\eta(\epsilon)>\bar{\eta}+\epsilon$. Thus for all $\epsilon<\epsilon'<1$ we have
\begin{equation} \label{eq.gltphidisc}
g(x)<\eta(\epsilon)\phi(\Gamma_{\eta(\epsilon)}^{-1}(\Gamma_{\eta(\epsilon)}(x)+\tau)), \quad x>x_2(\epsilon).
\end{equation}
By (\ref{eq.fsmallg2}) for every $\epsilon\in (0,1)$ there exists an $x_1(\epsilon)>0$ such that
\begin{equation} \label{eq.fltgdisc}
f(x)\leq \epsilon \eta(\epsilon) \phi(x), \quad x>x_1(\epsilon).
\end{equation}
%Define $y_2(\epsilon)=\Gamma_{\eta(\epsilon)}(x_2(\epsilon))+\tau$. Now, if $x>x_2(\epsilon)$ we have that $y:=\Gamma_{\eta(\epsilon)}(x)+\tau$
%obeys $y>y_2(\epsilon)$ and $x=\Gamma_{\eta(\epsilon)}^{-1}(y-\tau)$. Thus for $y>y_2(\epsilon)$, by (\ref{eq.gltphi}) we have
%$g(\Gamma_{\eta(\epsilon)}^{-1}(y-\tau))<\eta(\epsilon)\phi(\Gamma_{\eta(\epsilon)}^{-1}(y))$ for $y>y_2(\epsilon)$. Therefore
%\begin{equation} \label{eq.gltphi2}
%(1+\epsilon)g(\Gamma_{\eta(\epsilon)}^{-1}(x-\tau))<\eta(\epsilon)(1+\epsilon)\phi(\Gamma_{\eta(\epsilon)}^{-1}(x)), \quad x>y_2(\epsilon).
%\end{equation}
Define
\begin{equation} \label{def.cdisc}
c(\epsilon)=\Gamma_{\eta(\epsilon)}(\psi^\ast+x_1(\epsilon)+x_2(\epsilon))+(1+\epsilon)\tau,
\end{equation}
and define also
\begin{equation}\label{def.xepsupperdisc}
x_\epsilon(n)=\Gamma_{\eta(\epsilon)}^{-1}((1+\epsilon)nh+c(\epsilon)),\quad  n\geq -N.
\end{equation}
This function is well--defined since $c(\epsilon)>\Gamma_{\eta(\epsilon)}(\psi^\ast)+(1+\epsilon)\tau$, so $c(\epsilon)-(1+\epsilon)\tau>\Gamma_{\eta(\epsilon)}(\psi^\ast)$, or $x_{\epsilon}(n)>\psi^\ast$ for all $n\in \{-N,\ldots,0\}$. Since $c(\epsilon)>\Gamma_{\eta(\epsilon)}(x_1(\epsilon))+(1+\epsilon)\tau$ and $\Gamma_{\eta(\epsilon)}$ is increasing, $\Gamma_{\eta(\epsilon)}^{-1}(c(\epsilon)-(1+\epsilon)\tau)>x_1(\epsilon)$, so $x_\epsilon(n)>x_1(\epsilon)$ for all $n\geq -N$. Therefore by (\ref{eq.fltgdisc}), $f(x_\epsilon(n))\leq \epsilon\eta(\epsilon) \phi(x_\epsilon(n))$
for $n\geq 0$.
%g(\Gamma_{\eta(\varepsilon)}^{-1}(\Gamma_{\eta(\varepsilon)}(x_\epsilon(t))-\tau))
%=\epsilon g(\Gamma_{\eta(\epsilon)}^{-1}((1+\epsilon)t+c(\epsilon)-\tau))$.
Also for $n\geq 0$, we have
\begin{align*}
g(x_\epsilon(n-N))
&=g(\Gamma_{\eta(\epsilon)}^{-1}((1+\epsilon)(n-N)h+c(\epsilon))\\
&=g(\Gamma_{\eta(\epsilon)}^{-1}((1+\epsilon)nh-\tau-\epsilon \tau + c(\epsilon)))\\
&<g(\Gamma_{\eta(\epsilon)}^{-1}((1+\epsilon)nh-\tau+c(\epsilon))).
\end{align*}
Now, because $c(\epsilon)>\Gamma_{\eta(\epsilon)}(x_2(\epsilon))+\tau$, we have that the argument of $g$ on
the righthand side exceeds $x_2(\epsilon)$ for all $t\geq 0$. Therefore by \eqref{eq.gltphidisc}, we have
\begin{align*}
g(x_\epsilon(n-N))&<g(\Gamma_{\eta(\epsilon)}^{-1}((1+\epsilon)nh-\tau+c(\epsilon)))\\
&<\eta(\epsilon)\phi(
\Gamma_{\eta(\epsilon)}^{-1}(\Gamma_{\eta(\epsilon)}(\Gamma_{\eta(\epsilon)}^{-1}((1+\epsilon)nh-\tau+c(\epsilon)))+\tau)
)\\
&=\eta(\epsilon)\phi(\Gamma_{\eta(\epsilon)}^{-1}((1+\epsilon)nh-\tau+c(\epsilon))+\tau)\\
&=\eta(\epsilon)\phi(\Gamma_{\eta(\epsilon)}^{-1}((1+\epsilon)nh+c(\epsilon)))\\
&=\eta(\epsilon)\phi(x_\epsilon(n)).
\end{align*}
Hence
\begin{equation} \label{eq.fplusgdisc}
f(x_\epsilon(n))+g(x_\epsilon(n-N))<(1+\epsilon)\eta(\epsilon) \phi(x_\epsilon(n)),\quad n\geq 0.
\end{equation}
Now for $n\geq 0$, $\Gamma_{\eta(\epsilon)}(x_\epsilon(n))=(1+\epsilon)nh + c(\epsilon)$, so
\[
\Gamma_{\eta(\epsilon)}(x_\epsilon(n+1))-\Gamma_{\eta(\epsilon)}(x_\epsilon(n))=(1+\epsilon)h.
\]
Since $\Gamma_\eta$ is in $C^1$ and $(x_\epsilon(n))_{n\geq 0}$ is an increasing sequence, there exists
$\xi(n)\in [x_\epsilon(n),x_\epsilon(n+1)]$ such that
\[
\Gamma_{\eta(\epsilon)}(x_\epsilon(n+1))=\Gamma_{\eta(\epsilon)}(x_\epsilon(n))+\Gamma_{\eta(\epsilon)}'(\xi(n))
(x_\epsilon(n+1)-x_\epsilon(n)).
\]
Therefore we have
\[
(1+\epsilon)h=\Gamma_{\eta(\epsilon)}'(\xi(n)) (x_\epsilon(n+1)-x_\epsilon(n))=\frac{1}{\eta(\epsilon)}\frac{1}{\phi(\xi(n))}(x_\epsilon(n+1)-x_\epsilon(n)).
\]
Thus as $\phi$ is non--decreasing, as $\xi(n)\geq x_\epsilon(n)$, we have
\begin{equation}     \label{eq.xepsdiffereq1}
x_\epsilon(n+1)=x_\epsilon(n)+(1+\epsilon)\eta(\epsilon)h\phi(\xi(n))\geq x_\epsilon(n)+(1+\epsilon)\eta(\epsilon)h\phi(x_\epsilon(n)), \quad n\geq 0.
\end{equation}
Thus by \eqref{eq.fplusgdisc} and \eqref{eq.xepsdiffereq1} for $n\geq 0$ we have
\[
x_\epsilon(n+1)\geq x_\epsilon(n)+(1+\epsilon)\eta(\epsilon)h\phi(x_\epsilon(n))
> x_\epsilon(n) +  hf(x_\epsilon(n))+hg(x_\epsilon(n-N)).
\]
%\eta(\epsilon)(1+\epsilon)\phi(\Gamma_{\eta(\epsilon)}^{-1}((1+\epsilon)t+c(\epsilon))) =x_\epsilon'(t)$.
Now as $x_\epsilon(n)>\psi^\ast=\max_{n\in \{-N,\ldots,0\}}\psi(nh)$, we have $x_\epsilon(n)>x_h(n)$ for $n\in \{N,\ldots,0\}$.

Suppose that there is a $n_0\geq 1$ such that $x_\epsilon(n)>x_h(n)$ for $t\in \{-N,\ldots,n_0-1\}$
$x_\epsilon(n_0)\leq x_h(n_0)$. Therefore $x_\epsilon(n_0)-x_\epsilon(n_0-1)\leq x_h(n_0)-x_h(n_0-1)$.
Since $f$ and $g$ are non--decreasing,
\begin{align*}
x_\epsilon(n_0)-x_\epsilon(n_0-1)&\leq x_h(n_0)-x_h(n_0-1)\\
&=hf(x_h(n_0-1))+hg(x_h(n_0-1-N))\\
&\leq hf(x_\epsilon(n_0-1))+hg(x_h(n_0-N))\\
&\leq hf(x_\epsilon(n_0-1))+hg(x_\epsilon(n_0-1-N))\\
&<x_\epsilon(n_0)-x_\epsilon(n_0-1),
\end{align*}
a contradiction.

Thus $x_\epsilon(n)>x_h(n)$ for all $n\geq -N$. Hence $\Gamma_{\eta(\epsilon)}(x_h(n))<\Gamma_{\eta(\epsilon)}(x_\epsilon(n))$
for all $n\geq -N$. Hence
\[
\Gamma_{\eta(\epsilon)}(x_h(n))<\Gamma_{\eta(\epsilon)}(x_\epsilon(n))=(1+\epsilon)nh + c(\epsilon), \quad n \geq -N.
\]
But
$\Gamma(x_h(n))=\eta(\epsilon)\Gamma_{\eta(\epsilon)}(x_h(n))<(1+\epsilon)\eta(\epsilon)nh+\eta(\epsilon)c(\epsilon)$. Therefore
\[
\limsup_{n\rightarrow\infty} \frac{\Gamma(x_h(n))}{nh}\leq (1+\epsilon)\eta(\epsilon).
\]
Since $\epsilon>0$ is arbitrary, and $\eta(\epsilon)\rightarrow\eta$ as $\epsilon\rightarrow 0$, we have
\eqref{eq.growthlimsupdisc}. 

\subsection{Proof of Theorem~\ref{thm:liminfdisc}}
Suppose first that $\bar{\mu}_\epsilon$ is finite. Then
by (\ref{eq.gphigamma2}) for every $\epsilon\in (0,1)$ there exists $x_3(\epsilon)>0$ such that for $x>x_3(\epsilon)$
\begin{align*}
g(x)&>\bar{\mu}_\epsilon(1-\epsilon)\phi(\Gamma_{\mu(\epsilon)}^{-1}(\Gamma_{\mu(\epsilon)}(x)+(\tau+h)(1-\epsilon)))\\
&\geq \bar{\mu}(1-\epsilon)\phi(\Gamma_{\mu(\epsilon)}^{-1}(\Gamma_{\mu(\epsilon)}(x)+(\tau+h)(1-\epsilon)))\\
&>\mu(\epsilon)\phi(\Gamma_{\mu(\epsilon)}^{-1}(\Gamma_{\mu(\epsilon)}(x)+(\tau+h)(1-\epsilon))),
\end{align*}
where the penultimate inequality is a consequence of (\ref{eq.narmumu}), and the last inequality holds for all $\epsilon<\epsilon'$, because
for such $\epsilon$ we have $\mu(\epsilon)<(1-\epsilon)\bar{\mu}$. This holds for the following reason.

By \eqref{eq.muepstomu}, there exists $\epsilon_1\in (0,1)$ such that $\epsilon\in (0,\epsilon_1)$ implies
$-\epsilon<\mu(\epsilon)-\mu<\mu \epsilon$. Since $\mu<\bar{\mu}$, it follows that there exists $\epsilon_2\in (0,1)$ such that $\epsilon<\epsilon_2$ implies $\bar{\mu}>(1+\epsilon)\mu/(1-\epsilon)$. Hence for all $\epsilon<\epsilon':=\epsilon_1\wedge\epsilon_2$, we have
$\mu(\epsilon)<\mu(1+\epsilon)<(1-\epsilon)\bar{\mu}$.

Thus for all $0<\epsilon<\epsilon'<1$, and $x>x_3(\epsilon)$ we have
\begin{equation} \label{eq.ggtphidisc}
g(x)>\mu(\epsilon)\phi(\Gamma_{\mu(\epsilon)}^{-1}(\Gamma_{\mu(\epsilon)}(x)+(\tau+h)(1-\epsilon))), \quad x>x_3(\epsilon).
\end{equation}
When $\bar{\mu}_\epsilon=+\infty$, because $\mu(\epsilon)$ is finite, \eqref{eq.ggtphidisc}  is trivial.

Define $y_3(\epsilon)=\Gamma_{\mu(\epsilon)}(x_3(\epsilon))+(\tau+h)(1-\epsilon)$. Then for $y>y_3(\epsilon)$, if we define
$x=\Gamma_{\mu(\epsilon)}^{-1}(y-(\tau+h)(1-\epsilon))$, for $x>x_3(\epsilon)$ we have that $y>y_3(\epsilon)$. Thus by \eqref{eq.ggtphidisc}
\begin{equation}    \label{eq.ggtphi2disc}
g(\Gamma_{\mu(\epsilon)}^{-1}(y-(\tau+h)(1-\epsilon)))>\mu(\epsilon)\phi(\Gamma_{\mu(\epsilon)}^{-1}(y)), \quad y>y_3(\epsilon).
\end{equation}
Next let $N_0(\epsilon)=\inf\{n>0\,:\, x_h(n)\geq x_3(\epsilon)\}$ and define $N_1>N_0$ such that
\[
(1-\epsilon)(\tau+h)\Gamma_{\mu(\epsilon)}(x_h(N_0))\leq \Gamma_{\mu(\epsilon)}(x_h(N_1)).
\]
%, or $(1-\epsilon)\tau+  \Gamma_{\mu(\epsilon)}(x_3(\epsilon)))=\Gamma_{\mu(\epsilon)}(x(T_1(\epsilon)))$.
Define
\begin{equation}\label{def.xepslowerdisc}
x_\epsilon(n)=\Gamma_{\mu(\epsilon)}^{-1}((1-\epsilon)(n-N_1)h+\Gamma_{\mu(\epsilon)}(x_h(N_0))),\quad  n\geq N_1.
\end{equation}
Therefore for $n\geq N_1+N$ we have
\begin{align*}
(1-\epsilon)(n+1-N_1(\epsilon))h+\Gamma_{\mu(\epsilon)}(x_h(N_0))
&\geq (1-\epsilon)(\tau+h)+\Gamma_{\mu(\epsilon)}(x_h(N_0))\\
%=\Gamma_{\mu(\epsilon)}(x(T_1(\epsilon)))\\
&\geq (1-\epsilon)(\tau+h)+  \Gamma_{\mu(\epsilon)}(x_3(\epsilon)))
=y_3(\epsilon).
\end{align*}
Setting $y=(1-\epsilon)(n+1-N_1)h+\Gamma_{\mu(\epsilon)}(x_h(N_0))$ in
\eqref{eq.ggtphi2disc} yields
\begin{multline*}
g(\Gamma_{\mu(\epsilon)}^{-1}((1-\epsilon)(n-N_1-N)h+\Gamma_{\mu(\epsilon)}(x_h(N_0))))
\\>\mu(\epsilon)
\phi(\Gamma_{\mu(\epsilon)}^{-1}((1-\epsilon)(n+1-N_1)h+\Gamma_{\mu(\epsilon)}(x_h(N_0)))), \quad n\geq N_1+N.
\end{multline*}
By \eqref{def.xepslowerdisc} we have
\begin{equation} \label{eq.glowerdisc}
g(x_\epsilon(n-N))>\mu(\epsilon)\phi(x_\epsilon(n+1)), \quad n\geq N_1+N.
\end{equation}
Therefore by (\ref{eq.glowerdisc}) for $n\geq N_1+N$, and the fact that 
\[
\Gamma_{\mu(\epsilon)}(x_\epsilon(n))=(1-\epsilon)(n-N_1)h + \Gamma_{\mu(\epsilon)}(x_h(N_0)),
\] 
we have
\[
\Gamma_{\mu(\epsilon)}(x_\epsilon(n+1))-\Gamma_{\mu(\epsilon)}(x_\epsilon(n))
=(1-\epsilon)h.
\]
Hence there is $\xi(n)\in [x_\epsilon(n),x_\epsilon(n+1)]$ such that 
\[
x_\epsilon(n+1)-x_\epsilon(n)=(1-\epsilon)h \mu(\epsilon)\phi(\xi(n)).
\]
%$\Gamma_{\eta(\epsilon)}'(x_\epsilon(t)) x_\epsilon'(t)=(1+\epsilon)$, or
%$x_\epsilon'(t)=(1+\epsilon)\eta(\epsilon)\phi(x_\epsilon(t))$. Hence
Since $\phi$ is non--decreasing and $\xi(n)\leq x_\epsilon(n+1)$, we have 
\begin{align*} % \label{eq.xepsdiffeq2}
x_\epsilon(n+1)&=x_\epsilon(n)+(1-\epsilon)h \mu(\epsilon)\phi(\xi(n))\\
&\leq x_\epsilon(n)+(1-\epsilon)h \mu(\epsilon)\phi(x_\epsilon(n+1)).
\end{align*}
Therefore by \eqref{eq.glowerdisc}, we get for $n\geq N_1+N$
\begin{align*} 
x_\epsilon(n+1)&\leq x_\epsilon(n)+(1-\epsilon)h \mu(\epsilon)\phi(x_\epsilon(n+1))\\
&< x_\epsilon(n)+h(1-\epsilon)g(x_\epsilon(n-N))\\
&\leq x_\epsilon(n)+hf(x_\epsilon(n))+h(1-\epsilon)g(x_\epsilon(n-N))\\
&< x_\epsilon(n)+hf(x_\epsilon(n))+hg(x_\epsilon(n-N)).
\end{align*}
Now for $n\in \{N_1,\ldots, N_1+N\}$ we have
\begin{align*}
x_\epsilon(n)
&\leq x_\epsilon(N_1+N)
=\Gamma_{\mu(\epsilon)}^{-1}((1-\epsilon)\tau+\Gamma_{\mu(\epsilon)}(x_h(N_0)))\\
&<\Gamma_{\mu(\epsilon)}^{-1}((1-\epsilon)(\tau+h)+\Gamma_{\mu(\epsilon)}(x_h(N_0)))\\
&\leq \Gamma_{\mu(\epsilon)}^{-1}(\Gamma_{\mu(\epsilon)}(x_h(N_1)))=x_h(N_1)\leq x_h(n),
\end{align*}
where we used at the last step the fact that $x_h$ is increasing on $\{N_1,\ldots,N_1+N\}\subset\{N,N+1,\ldots\}$.
%Finally $x_\epsilon(T_1(\epsilon))<x_\epsilon(T_1(\epsilon)+\tau)=x(T_1(\epsilon))$.
Therefore we have $x_\epsilon(n)<x_h(n)$ for $n\in \{N_1(\epsilon),\ldots, N_1(\epsilon)+N\}$, and also  
$x_\epsilon(n+1)< x_\epsilon(n)+hf(x_\epsilon(n))+hg(x_\epsilon(n-N))$ for $n\geq N_1+N$.
%<f(x_\epsilon(t))+g(x_\epsilon(t-\tau))$ for $t\geq T_1+\tau$.

Suppose that there is a $n_1\geq N_1(\epsilon)+N+1$ such that $x_\epsilon(n)<x_h(n)$ for $n\in \{N_1(\epsilon),\ldots,n_1\}$ and $x_\epsilon(n_1)\geq x_h(n_1)$. Therefore $x_\epsilon(n_1)-x_\epsilon(n_1-1)\geq x_h(n_1)-x_h(n_1-1)$. Then as
$f$ and $g$ are non--decreasing,
\begin{align*}
x_\epsilon(n_1)-x_\epsilon(n_1-1) &\geq x_h(n_1)-x_h(n_1-1)\\
&=hf(x_h(n_1-1))+hg(x_h(n_1-1-N))\\
%x_\epsilon'(t_1)&\geq x'(t_1)=f(x(t_1))+g(x(t_1-\tau))\\
%&=f(x_\epsilon(t_1))+g(x(t_1-\tau))
&\geq hf(x_\epsilon(n_1-1))+hg(x_\epsilon(n_1-1-N))\\
&>x_\epsilon(n_1) -x_\epsilon(n_1-1),
%&\geq f(x_\epsilon(t_1))+g(x_\epsilon(t_1-\tau))>x_\epsilon'(t_1),
\end{align*}
a contradiction. Thus $x_\epsilon(n)<x_h(n)$ for all $n\geq N_1$. Hence
$\Gamma_{\mu(\epsilon)}(x_h(n))>\Gamma_{\mu(\epsilon)}(x_\epsilon(n))$
for all $n\geq N_1(\epsilon)$. Hence
\[
\Gamma_{\mu(\epsilon)}(x_h(n))>\Gamma_{\mu(\epsilon)}(x_\epsilon(n))
=(1-\epsilon)(n-N_1)+ \Gamma_{\mu(\epsilon)}(x_h(N_0)), \quad n \geq N_1(\epsilon).
\]
But $\Gamma(x_h(n))
=\mu(\epsilon)\Gamma_{\mu(\epsilon)}(x_h(n))
>(1-\epsilon)\mu(\epsilon)n+\mu(\epsilon)\Gamma_{\mu(\epsilon)}(x_h(N_0))$. 
Therefore
\[
\liminf_{n\rightarrow\infty} \frac{\Gamma(x_h(n))}{nh} \geq (1-\epsilon)\mu(\epsilon).
\] 
Since $\epsilon>0$ is arbitrary, and $\mu(\epsilon)\rightarrow\mu$ as $\epsilon\rightarrow 0$, we have (\ref{eq.growthliminfdisc}).

\subsection{Proof of Theorem~\ref{thm:presrv1}}
Let $j\geq N$. Summing across both sides of \eqref{eq:disceq} yields
\[
x_h(j+1)=x_h(j-N)+h\sum_{n=j-N}^{j} f(x_h(n))+h\sum_{n=j-N}^{j}g(x_h(n-N)).
\]
Let $\epsilon(\tau+h)<1/2$. Since $x_h(n)\to\infty$ as $n\to\infty$ and $f(x)/x\to 0$ as $x\to\infty$, there exists $N_1(\epsilon)$ such that $f(x_h(n))\leq \epsilon x_h(n)$ for all $n\geq N_1(\epsilon)$. Hence for $j\geq N_1(\epsilon)$ we have
\begin{align*}
x_h(j+1)&\leq x_h(j-N)+h\sum_{n=j-N}^{j} \epsilon x_h(n) +h\sum_{n=j-N}^{j}g(x_h(n-N))\\
&\leq x_h(j-N)+ h(N+1)\epsilon x_h(j) +h\sum_{n=j-N}^{j}g(x_h(n-N))\\
&\leq x_h(j-N)+ h(N+1)\epsilon x_h(j+1) +h\sum_{n=j-N}^{j}g(x_h(n-N)).
\end{align*}
Hence for $j\geq N_1(\epsilon)$ we have
\[
x_h(j+1) \leq \frac{1}{1-(\tau+h)\epsilon}x_h(j-N)+\frac{h}{1-(\tau+h)\epsilon}h\sum_{n=j-N}^{j}g(x_h(n-N)).
\]
Since $g$ is in $\text{RV}_\infty(1)$, $x\mapsto g(x)/x$ is asymptotic to a non--decreasing function,
there exists $g_0$ such that $g_0$ is non--decreasing, $g_0(x)\to \infty$ as $x\to\infty$
and $g_0(x)/g(x)/x\to 1$ as $x\to\infty$. Therefore $g_1$ defined by $g_1(x):=xg_0(x)$ is increasing and is
in $\text{RV}_\infty(1)$. Since $x_h(n)\to\infty$ as $n\to\infty$, for every $\epsilon>0$ there exists $N_2(\epsilon)\geq N$ such that $g(x_h(n-N))<(1+\epsilon)g_1(x_h(n-N))$. Thus for $j\geq N_2(\epsilon)$ we have
\begin{align*}
h\sum_{n=j-N}^{j}g(x_h(n-N))&\leq h(1+\epsilon)\sum_{n=j-N}^{j}g_1(x_h(n-N))\\
&\leq h(N+1)(1+\epsilon) g_1(x_h(j-N)).
\end{align*}
Hence
\[
h\sum_{n=j-N}^{j}g(x_h(n-N)) \leq (\tau+h)(1+\epsilon) g_1(x_h(j-N)), \quad j\geq N_2(\epsilon).
\]
Let $N_3=\max(N_1,N_2)$. Then for $j\geq N_3$ we have
\[
x_h(j+1) \leq x_h(j-N) + \left(\frac{1}{1-(\tau+h)\epsilon}-1\right)x_h(j-N)
+\frac{(\tau+h)(1+\epsilon)}{1-(\tau+h)\epsilon} g_1(x_h(j-N)).
\]
Define $x^\ast_h(n)=x_h(n(N+1))$ for $n\geq 0$. Therefore for $n\geq N_3$ we have
\[
x^\ast_h(j+1) \leq x^\ast_h(j) + \left(\frac{1}{1-(\tau+h)\epsilon}-1\right)x^\ast_h(j)+\frac{(\tau+h)(1+\epsilon)}{1-(\tau+h)\epsilon} g_1(x^\ast_h(j)).
\]
Define
\begin{equation}  \label{def.geps}
g_\epsilon(x)=
\left(\frac{1}{1-(\tau+h)\epsilon}-1\right)x+\frac{(\tau+h)(1+\epsilon)}{1-(\tau+h)\epsilon} g_1(x), \quad x> 0.
\end{equation}
Then $g_\epsilon$ is in $\text{RV}_\infty(1)$, $x\mapsto g_\epsilon(x)/x$ is positive and non--decreasing, and $g_\epsilon(x)/x\to\infty$ as $x\to\infty$. Moreover
\[
x^\ast_h(n+1) \leq x^\ast_h(n) + g_\epsilon(x^\ast_h(n)), \quad n\geq N_3(\epsilon).
\]
Next, define
\[
y_\epsilon(n+1) = y_\epsilon(n) + g_\epsilon(y_\epsilon(n)), \quad n\geq N_3(\epsilon); \quad
y_\epsilon(N_3)=2x^\ast_h(N_3(\epsilon)).
\]
Since $g_\epsilon$ is increasing, it follows that $x^\ast_h(n)\leq y_\epsilon(n)$ for all $n\geq N_3(\epsilon)$.
Define
\[
H_\epsilon(x)=\int_1^x \frac{1}{u\log(1+g_\epsilon(u)/u)}\,du, \quad x\geq 0.
\]
Then by applying Lemma~\ref{lemma.l1} to $(y_\epsilon)$, we have that
\[
\lim_{n\to\infty} \frac{H_\epsilon(y_\epsilon(n))}{n}=1.
\]
Since $H_\epsilon$ is increasing, and $x^\ast_h(n)\leq y_\epsilon(n)$ for all $n\geq N_3(\epsilon)$, we have by
the definition of $x^\ast_h$ that
\[
\limsup_{n\to\infty} \frac{H_\epsilon(x_h(n(N+1)))}{n}
=
\limsup_{n\to\infty} \frac{H_\epsilon(x^\ast_h(n))}{n}
\leq   \lim_{n\to\infty} \frac{H_\epsilon(y_\epsilon(n))}{n}=1.
\]
Now by L'H\^opital's rule and \eqref{def.geps}
\[
\lim_{x\to\infty} \frac{H_\epsilon(x)}{G(x)}
=
\lim_{x\to\infty} \frac{\log(1+g(x)/x)}{\log(1+g_\epsilon(x)/x)}
=
\lim_{x\to\infty} \frac{\log(1+g(x)/x)}{\log(\frac{1}{1-(\tau+h)\epsilon}+\frac{(\tau+h)(1+\epsilon)}{1-(\tau+h)\epsilon} \frac{g_1(x)}{x})}.
\]
Since $g(x)/g_1(x)\to 1$ as $x\to\infty$, we have that
\[
\lim_{x\to\infty} \frac{\log(1+g(x)/x)}{\log(1+g_1(x)/x)}=1.
\]
Therefore $\lim_{x\to\infty} H_\epsilon(x)/G(x)=1$. Hence
\begin{equation} \label{eq.GxNn}
\limsup_{n\to\infty} \frac{G(x_h(n(N+1)))}{n}\leq 1.
\end{equation}
For every $n\in \mathbb{N}$ there exists $j=j(n)\geq 1$ such that $n(N+1)\leq j<(n+1)(N+1)$. Since $G$ is increasing,
and $(x_h(n))_{n\geq 0}$ is increasing, we have
\begin{align*}
\frac{G(x_h(j))}{jh}&\leq \frac{G(x_h((n+1)(N+1)))}{jh} \\
&\leq \frac{G(x_h((n+1)(N+1)))}{n(N+1)h}\\
&=\frac{1}{\tau+h}\frac{G(x_h((n+1)(N+1)))}{n+1}\cdot\frac{n+1}{n}.
\end{align*}
By \eqref{eq.GxNn}, we have
\[
\limsup_{j\to\infty} \frac{G(x_h(j))}{jh} \leq \frac{1}{\tau+h},
\]
which gives the desired upper limit in \eqref{eq.growthlimdelaydisc}.

To get a lower bound, since $f(x)\geq 0$, we have
$x_h(n+1)\geq x_h(n)+hg(x_h(n-N))$ for $n\geq 0$. Since $x_h(n)\to\infty$ as
$n\to\infty$, for every $\epsilon\in (0,1)$ there exists $N_4(\epsilon)\geq N$ such that $g(x_h(n-N))>(1-\epsilon)g_1(x_h(n-N))$. Let $N_5(\epsilon)=\max(N_4(\epsilon),N_1(\epsilon))$. Let $y_h^{(1)}$ be defined by
\begin{align*}
y_h^{(1)}(n+1)&= y_h^{(1)}(n)+h(1-\epsilon)g_1(y_h^{(1)}(n-N)), \quad n\geq N_5(\epsilon); \\
y_h^{(1)}(n)&=x_h(n)/2, \quad n=N_5(\epsilon)-N,\ldots,N_5(\epsilon).
\end{align*}
Then we have for $n\geq N_5(\epsilon)$ the inequality $x_h(n+1)\geq x_h(n)+h(1-\epsilon)g_1(x_h(n-N))$. Hence
$y_h^{(1)}(n)\leq x_h(n)$ for $n\geq N_5(\epsilon)-N$. Clearly $(y_h^{(1)}(n))_{n\geq N_5(\epsilon)}$ is increasing
and $y_h^{(1)}(n)\to\infty$ as $n\to\infty$.

Let $n\geq N_5(\epsilon)+N$. Then as $y_h^{(1)}$ is increasing, we have
\[
y_h^{(1)}(n+1)= y_h^{(1)}(n)+h(1-\epsilon)g_1(y_h^{(1)}(n-N))
\geq  y_h^{(1)}(n-N)+h(1-\epsilon)g_1(y_h^{(1)}(n-N)).
\]
Therefore for $n\geq N_5(\epsilon)+N$ we have
\begin{multline*}
\log y_h^{(1)}(n+1) \geq \log \left(\frac{g_1(y_h^{(1)}(n-N))}{y_h^{(1)}(n-N)}\right)+\log y_h^{(1)}(n-N)\\+
\log\left( h(1-\epsilon)+ \frac{y_h^{(1)}(n-N)}{g_1(y_h^{(1)}(n-N))}\right),
\end{multline*}
and so
\begin{equation*}
\log y_h^{(1)}(n+1) \geq \log y_h^{(1)}(n-N)+\log( h(1-\epsilon))+\log g_0(y_h^{(1)}(n-N)).
\end{equation*}
Define $u(n):=\log y_h^{(1)}(n)$ for $n\geq N_5(\epsilon)$. Then $(u(n))_{n\geq N_5}$ is increasing and tends to infinity as $n\to\infty$, and with $\gamma_0(x):=\log( h(1-\epsilon))+\log g_0(e^{x})$, we have
\[
u(n+1) \geq u(n-N)+\gamma_0(u(n-N)), \quad n\geq N_5(\epsilon)+N.
\]
Since $g_0$ is non--decreasing, so is $\gamma_0$, and moreover $\gamma_0(x)\to \infty$ as $x\to\infty$.
Since $g_0$ is in $\text{RV}_\infty(0)$, there is $g_3$ in $\text{RV}_\infty(0)$ which is also in $C^1$ such that $g(x)/g_3(x)\to 1$ as $x\to\infty$, $xg_3'(x)/g_3(x)\to 0$ as $x\to\infty$. Clearly for $x^\ast$ sufficiently large we have $g_3(e^{x})>e$ for all $x>x^\ast$, and so we may define
\[
G_3(x)=\int_{x^\ast}^x \frac{1}{\log g_3(e^{u})}\,du.
\]
Then $G_3'(x)=1/\log g_3(e^{x})>0$ for $x>x^\ast$ and since $g_3$ is in $C^1$ we have
\[
G_3''(x)=-\frac{d}{dx} \log g_3(e^{x}) \cdot \frac{1}{(\log g_3(e^{x}))^2}
=-\frac{1}{g_3(e^{x}))}g_3'(e^x) e^x \cdot \frac{1}{(\log g_3(e^{x}))^2}.
\]

Since there $u(n)\to \infty$, there is $N_6$ is such that $u(n)>x^\ast$ for $n\geq N_6$. Let $N_7(\epsilon)=\max(N_5(\epsilon),N_6)+N$.
Then for $n\geq N_7(\epsilon)$ we have $G_3(u(n+1)) \geq G_3(u(n-N)+\gamma_0(u(n-N)))$ and so by Taylor's theorem, there exists $\xi_n\in [u(n-N),u(n-N)+\gamma_0(u(n-N))]$ such that
\begin{align*}
\lefteqn{G_3(u(n+1))}\\
&\geq G_3(u(n-N)+\gamma_0(u(n-N)))\\
&=G_3(u(n-N))+G_3'(u(n-N))\gamma_0(u(n-N))+\frac{1}{2}G_3''(\xi_n)\gamma_0^2(u(n-N)),
\end{align*}
for $n\geq N_7(\epsilon)$. Next, with $\eta_n:=g_3'(e^\xi_n) e^{\xi_n}/g_3(e^{\xi_n}))$ and using the fact that $xg_3'(x)/g_3(x)\to 0$ as $x\to\infty$, we have that $\eta_n\to 0$ as $n\to\infty$. Define for $n\geq N_7(\epsilon)$ the sequence
\begin{multline*}
\delta(n):=\frac{\log( h(1-\epsilon))+\log g_0(e^{u(n-N)})}{\log g_3(e^{u(n-N)})} -1\\
-\frac{1}{2} \eta_n \frac{\left(\log( h(1-\epsilon))+\log g_0(e^{u(n-N)})\right)^2}{(\log g_3(e^{\xi_n}))^2}.
\end{multline*}
so that
\[
G_3(u(n+1)) \geq G_3(u(n-N))+1+\delta(n), \quad n\geq N_7(\epsilon).
\]

Since $\xi_n\to\infty$ as $n\to\infty$ and $g_3(x)/g_0(x)\to 1$ as $x\to\infty$ we have that  for every
$\epsilon\in (0,1)$ that there exists $N_8(\epsilon)$ such that
$\log g_3(e^{\xi_n}) >  \log(1-\epsilon)+\log g_0(e^{\xi_n})$ for all $n\geq N_8(\epsilon)$ and so for $n\geq N_9(\epsilon)=\max(N_8(\epsilon),N_7(\epsilon))+N$
and so
\[
\frac{\left(\log( h(1-\epsilon))+\log g_0(e^{u(n-N)})\right)^2}{(\log g_3(e^{\xi_n}))^2}
\leq \frac{\left(\log( h(1-\epsilon))+\log g_0(e^{u(n-N)})\right)^2}{(\log(1-\epsilon)+\log g_0(e^{\xi_n}))^2}.
\]
Since $g_0$ is increasing and $\xi_n\geq u(n-N)$ we have $\log g_0(e^{\xi_n})\geq \log g_0(e^{u(n-N)})$.
Hence
\[
\frac{\left(\log( h(1-\epsilon))+\log g_0(e^{u(n-N)})\right)^2}{(\log g_3(e^{\xi_n}))^2}
\leq \frac{\left(\log( h(1-\epsilon))+\log g_0(e^{\xi_n})\right)^2}{(\log(1-\epsilon)+\log g_0(e^{\xi_n}))^2}.
\]
Therefore
\[
\limsup_{n\to\infty}
\frac{\left(\log( h(1-\epsilon))+\log g_0(e^{u(n-N)})\right)^2}{(\log g_3(e^{\xi_n}))^2}\leq 1,
\]
 and so $\delta(n)\to 0$ as $n\to\infty$. Let $z(n)=G_3(u(n))$. Note that $z$ is increasing and $z(n)\to\infty$ as $n\to\infty$. Then we have $z(n+1) \geq z(n-N)+1+\delta(n)$. Let $j\in \{0,\ldots,N\}$.
Define $z_j^\ast(n)=z((N+1)n+j)$. Then
\begin{align*}
z_j^\ast(n)&=z(Nn+n+j-1+1)\\
&\geq z(Nn+n+j-1-N)+1+\delta(Nn+n+j-1)\\
&=z_j^\ast(n-1)+1+\delta(Nn+n+j-1).
\end{align*}
Now for $n\geq n'$ we have
\[
\sum_{m=n'}^n z_j^\ast(m) \geq \sum_{m=n'}^n z_j^\ast(m-1)+n-n'+1+\sum_{m=n'}^{n}\delta(Nm+m+j-1),
\]
so
\[
\frac{z_j^\ast(n)}{n} \geq \frac{z_j^\ast(n'-1)}{n}+1 + \frac{-n'+1}{n}+\frac{1}{n}\sum_{m=n'}^{n}\delta(Nm+m+j-1).
\]
Since $\delta(n)\to 0$ as $n\to\infty$, we have $\liminf_{n\to\infty} z_j^\ast(n)/n\geq 1$. Therefore
\[
\liminf_{n\to\infty} \frac{z((N+1)n+j)}{n(N+1)}\geq \frac{1}{N+1}, \quad\text{for each $j=0,\ldots,N$}.
\]
Hence
\[
\liminf_{n\to\infty} \frac{G_3(\log y_h^{(1)}(n))}{n}
=
\liminf_{n\to\infty} \frac{G_3(u(n))}{n}
=
\liminf_{n\to\infty} \frac{z(n)}{n}\geq \frac{1}{N+1}.
\]
Since $x_h(n)\geq y_h^{(1)}(n)$ for $n\geq N_5(\epsilon)-N$, and $G_3$ is increasing, we have
\begin{equation} \label{eq.Glogxnhlower}
\liminf_{n\to\infty} \frac{G_3(\log x_h(n))}{nh}\geq
\liminf_{n\to\infty} \frac{G_3(\log y_h^{(1)}(n))}{nh}\geq \frac{1}{Nh+h}=\frac{1}{\tau+h}.
\end{equation}
Now
\begin{equation} \label{eq.defG4}
G_3(\log x))=\int_{x^\ast}^{\log x}  \frac{1}{\log g_3(e^v)}\,dv
=\int_{e^{x^\ast}}^{x}  \frac{1}{u\log g_3(u)}\,du=:G_4(x).   % u=e^v => log u=v
\end{equation}
Since $g_3(x)/g_0(x)\to 1$ as $x\to\infty$ and each belongs to $\text{RV}_\infty(0)$, we have that
\[
\lim_{x\to\infty} \frac{\log g_0(x)}{\log g_3(x)}=1.
\]
Similarly, as $(1+g(x)/x)/g_0(x)\to 1$ as $x\to\infty$ and $g_0$ is in $\text{RV}_\infty(0)$,
\[
\lim_{x\to\infty} \frac{\log (1+g(x)/x)}{\log g_0(x)}=1.
\]
Using these limits and L'H\^opital's rule, we arrive at
\[
\lim_{x\to\infty} \frac{G_4(x)}{G(x)}
=\lim_{x\to\infty} \frac{\log(1+g(x)/x)}{\log g_3(x)}
=\lim_{x\to\infty} \frac{\log(1+g(x)/x)}{\log g_0(x)}\cdot \frac{\log g_0(x)}{\log g_3(x)}=1.
\]
Since $x_h(n)\to\infty$ as $n\to\infty$ and \eqref{eq.Glogxnhlower} and $G_4$ is defined by \eqref{eq.defG4},
by using the last limit, we get
\[
\liminf_{n\to\infty} \frac{G(x_h(n))}{nh}
=\liminf_{n\to\infty} \frac{G(x_h(n))}{G_4(x_h(n))}\frac{G_4(x_h(n))}{nh}
=\liminf_{n\to\infty} \frac{G_3(\log x_h(n))}{nh}\geq \frac{1}{\tau+h},
\]
which is the lower limit in \eqref{eq.growthlimdelaydisc}.

In order to prove \eqref{eq.growthlimdelaydiscrv1interp}, notice for any $t>0$
that there exists $n\geq 0$ such that $nh\leq t<(n+1)h$. Also as the linear interpolant $\bar{x}_h$ defined
by \eqref{def.lininterpolant}, we have $x_h(n)\leq \bar{x}_h(t)\leq x_h(n+1)$. Therefore
\[
\frac{G(\bar{x}_h(t))}{t}\leq \frac{G(x_h(n+1))}{nh}=\frac{G(x_h(n+1))}{(n+1)h}\cdot\frac{n+1}{n}.
\]
Therefore by \eqref{eq.growthlimdelaydisc}, we have
\begin{equation} \label{eq.limsupinterprv1}
\limsup_{t\to\infty} \frac{G(\bar{x}_h(t))}{t}\leq \frac{1}{\tau+h}.
\end{equation}
To get the lower bound, we observe that for $nh\leq t<(n+1)h$, we have
\[
\frac{G(\bar{x}_h(t))}{t}\geq \frac{G(x_h(n))}{(n+1)h}=\frac{G(x_h(n))}{nh}\cdot \frac{n}{n+1}.
\]
Therefore by \eqref{eq.growthlimdelaydisc}, we have
\[
\liminf_{t\to\infty} \frac{G(\bar{x}_h(t))}{t}\geq \frac{1}{\tau+h}.
\]
Combining this limit with \eqref{eq.limsupinterprv1} yields \eqref{eq.growthlimdelaydiscrv1interp}.

\subsection{Proof of Theorem~\ref{thm:limnoinst}}
Let $N \in\mathbb{N}$ and set $h=\tau/N$. Let $j\geq N$. Integrating
over $[(j-N)h,(j+1)h]$ yields
\[
x((j+1)h)=x((j-N)h)+\int_{(j-N)h}^{(j+1)h} f(x(s))\,ds
+\int_{(j-N)h}^{(j+1)h} g(x(s-Nh))\,ds.
\]
Let $\epsilon(\tau+ h)<1/2$. Since $x(t)\to\infty$ as $t\to\infty$
and $f(x)/x\to 0$ as $x\to\infty$, there exists $T_1(\epsilon)>\tau$
such that $f(x(s))\leq \epsilon x(s)$ for all $s\geq T_1(\epsilon)$.
Let $N_1(\epsilon)$ be an integer such that $N_1(\epsilon)
h>T_1(\epsilon)$. Then for $j\geq N_1(\epsilon)$, and using the fact
that $x$ is increasing, we have
\begin{align*}
x((j+1)h)&\leq x((j-N)h)+\int_{(j-N)h}^{(j+1)h} \epsilon x(s)\,ds + \int_{(j-N)h}^{(j+1)h} g(x(s-Nh))\,ds\\
&\leq x((j-N)h)+h(N+1)\epsilon x((j+1)h) +\int_{(j-N)h}^{(j+1)h}
g(x(s-Nh))\,ds\\
&\leq x((j-N)h)+(\tau+h)\epsilon x((j+1)h) +\int_{(j-N)h}^{(j+1)h}
g(x(s-Nh))\,ds.
\end{align*}
Hence for $j\geq N_1(\epsilon)$ we have
\begin{multline*}
x((j+1)h) \leq x((j-N)h) +
\left(\frac{1}{1-(\tau+h)\epsilon}-1\right)x((j-N)h)\\
+\frac{1}{1-(h+\tau)\epsilon} \int_{(j-N)h}^{(j+1)h} g(x(s-Nh))\,ds.
\end{multline*}
Since $g$ is in $\text{RV}_\infty(1)$, $x\mapsto g(x)/x$ is
asymptotic to a non--decreasing function, there exists $g_0$ such
that $g_0$ is non--decreasing, $g_0(x)\to \infty$ as $x\to\infty$
and $g_0(x)/g(x)/x\to 1$ as $x\to\infty$. Therefore $g_1$ defined by
$g_1(x):=xg_0(x)$ is increasing and is in $\text{RV}_\infty(1)$.
Since $x(t)\to\infty$ as $t\to\infty$, for every $\epsilon>0$ there
exists $T_2(\epsilon)\geq \tau$ such that
$g(x(t-\tau))<(1+\epsilon)g_1(x(t-\tau))$ for all $t\geq T_2(\epsilon)$.
Let $N_2(\epsilon)$ be an integer such that $N_2(\epsilon) h>T_2(\epsilon)$. Thus for $j\geq N_2(\epsilon)+N$ we have
$jh\geq N_2(\epsilon)h+Nh>T_2+\tau\geq 2\tau=2Nh$, so as $x$ is
increasing on $[0,\infty)$ we have
\begin{align*}
\int_{(j-N)h}^{(j+1)h} g(x(s-Nh))\,ds &\leq (1+\epsilon)\int_{(j-N)h}^{(j+1)h} g_1(x(s-Nh))\,ds\\
&\leq h(N+1)(1+\epsilon) g_1(x((j+1-N)h)).
\end{align*}
Let $N_3(\epsilon)=\max(N_1(\epsilon),N_2(\epsilon)+N)$. Then for
$j\geq N_3(\epsilon)$ we have
\begin{multline*}
x((j+1)h) \leq x((j-N)h) +
\left(\frac{1}{1-(h+\tau)\epsilon}-1\right)x((j-N)h)\\
+\frac{(h+\tau)(1+\epsilon)}{1-(h+\tau)\epsilon}g_1(x(j+1-N)h)),
\end{multline*}
which, as $x$ is increasing, implies
\begin{multline*}
x((j+1)h) \leq x((j+1-N)h) +
\left(\frac{1}{1-(h+\tau)\epsilon}-1\right)x((j+1-N)h)\\+\frac{(h+\tau)(1+\epsilon)}{1-(h+\tau)\epsilon}
g_1(x(j+1-N)h)), \quad j\geq N_3(\epsilon).
\end{multline*}
%Thus, with $x_n=x(nh)$ we have
%\[
%x_{j+1} \leq x_{j+1-N} +
%\left(\frac{1}{1-(h+\tau)\epsilon}-1\right)x_{j+1-N}+\frac{1}{1-(h+\tau)\epsilon}(h+\tau)(1+\epsilon)
%g_1(x_{j+1-N}), \quad j\geq N_3(\epsilon).
%\]
Define $x^\ast_h(n)=x(nNh)$ for $n\geq -1$. Therefore for $n\geq
N_3$, and since $N\geq 1$ we have
\[
x^\ast_h(j+1) \leq x^\ast_h(j) +
\left(\frac{1}{1-(h+\tau)\epsilon}-1\right)x^\ast_h(j)+\frac{(h+\tau)(1+\epsilon)}{1-(h+\tau)\epsilon}
g_1(x^\ast_h(j)).
\]
The proof now continues as in the proof of
Theorem~\ref{thm:presrv1}, where $\tau$ is replaced by $\tau+h$.
Proceeding in this manner we arrive at
\begin{equation} \label{eq.GxNnh}
\limsup_{n\to\infty} \frac{G(x(nNh))}{n}\leq 1.
\end{equation}
For every $t>0$ there exists $n\in \mathbb{N}$ such that $nNh\leq t<(n+1)Nh$. Since $G$ is increasing,
and $x$ is increasing, we have
\[
\frac{G(x(t))}{t}\leq \frac{G(x((n+1)Nh))}{t} \leq \frac{G(x((n+1)Nh))}{nNh}
=\frac{1}{\tau}\frac{G(x((n+1)Nh))}{n+1}\cdot\frac{n+1}{n}.
\]
By \eqref{eq.GxNnh}, we have
\[
\limsup_{t\to\infty} \frac{G(x(t))}{t} \leq \frac{1}{\tau},
\]
%Setting $N=1$ gives $h=\tau$,
and therefore the desired upper limit in \eqref{eq.growthlimdelay}.

To get a lower bound, since $f(x)\geq 0$, we have
\[
x((n+1)h)\geq x(nh)+\int_{nh}^{(n+1)h} g(x(s-Nh))\,ds, \quad n\geq 0.
\]
Since $x(t)\to\infty$ as $t\to\infty$, for every $\epsilon\in (0,1)$ there exists $T_4(\epsilon)\geq \tau$
such that $g(x(t-\tau))>(1-\epsilon)g_1(x(t-\tau))$. Let $N_4(\epsilon)$ be ahn integer such that $N_4(\epsilon) h>T_4(\epsilon)$.
Let $N_5(\epsilon)=\max(N_4(\epsilon),N_1(\epsilon))$.
Thus for $n\geq N_5(\epsilon)$ we have $nh\geq N_5(\epsilon))h\geq \max(T_4(\epsilon),\tau)$, so as $x$ is increasing on $[0,\infty)$
we have
\begin{align*}
x((n+1)h)&\geq x(nh)+(1-\epsilon)\int_{nh}^{(n+1)h} g_1(x(s-Nh))\,ds\\
&\geq x(nh)+(1-\epsilon)h  g_1(x(nh-Nh)).
\end{align*}
Then with $x_h(n):=x(nh)$, we have the inequality
\[
x_h(n+1)\geq x_h(n)+(1-\epsilon)h  g_1(x_h(n-N)), \quad n\geq N_5(\epsilon).
\]
Let $y_h^{(1)}$ be defined by
\begin{align*}
y_h^{(1)}(n+1)&= y_h^{(1)}(n)+h(1-\epsilon)g_1(y_h^{(1)}(n-N)), \quad n\geq N_5(\epsilon); \\
y_h^{(1)}(n)&=x(nh)/2, \quad n=N_5(\epsilon)-N,\ldots,N_5(\epsilon).
\end{align*}
Hence $y_h^{(1)}(n)\leq x(nh)$ for $n\geq N_5(\epsilon)-N$. The proof now proceeds exactly as in Theorem~\ref{thm:presrv1}, and we arrive at the analogue of \eqref{eq.Glogxnhlower}, namely
\begin{equation} \label{eq.Glogxnhlowerh}
\liminf_{n\to\infty} \frac{G_3(\log x(nh))}{nh}\geq
\liminf_{n\to\infty} \frac{G_3(\log y_h^{(1)}(n))}{nh}\geq \frac{1}{Nh+h}=\frac{1}{\tau+h},
\end{equation}
where we have used the fact that $x(nh)=x_h(n)$. By \eqref{eq.defG4}, we have $G_3(\log x)=G_4(x)$, so once again we have that $\lim_{x\to\infty}G_4(x)/G(x)=1$. Since $x(nh)\to\infty$ as $n\to\infty$, \eqref{eq.Glogxnhlowerh} holds, and $G_4$ is defined by \eqref{eq.defG4}, by using the last limit, we get
\[
\liminf_{n\to\infty} \frac{G(x(nh))}{nh}
=\liminf_{n\to\infty} \frac{G(x(nh))}{G_4(x(nh))}\frac{G_4(x(nh))}{nh}
=\liminf_{n\to\infty} \frac{G_3(\log x(nh))}{nh}\geq \frac{1}{\tau+h}.
\]
Now, for every $t>0$ there exists $n$ such that $nh\leq t<(n+1)h$. Since $x$ is increasing and $G$ is increasing, we
have
\[
\frac{G(x(t))}{t} \geq \frac{G(x(nh))}{t}\geq \frac{G(x(nh))}{(n+1)h}= \frac{G(x(nh))}{nh}\frac{n}{n+1}.
\]
Therefore
\[
\liminf_{t\to\infty} \frac{G(x(t))}{t}\geq \liminf_{n\to\infty}  \frac{G(x(nh))}{nh}\geq \frac{1}{\tau+h}.
\]
Letting $h\to 0$ yields
\[
\liminf_{t\to\infty} \frac{G(x(t))}{t}\geq \frac{1}{\tau},
\]
which is the lower limit in \eqref{eq.growthlimdelay}.

\subsection{Proof of Theorem~\ref{thm:presrvbeta}}
Let $j\geq N$. Summing across both sides of \eqref{eq:disceq} yields
\[
x_h(j+1)=x_h(j-N)+h\sum_{n=j-N}^{j} f(x_h(n))+h\sum_{n=j-N}^{j}g(x_h(n-N)).
\]
Let $\epsilon(\tau+h)<1/2$. Since $x_h(n)\to\infty$ as $n\to\infty$ and $f(x)/x\to 0$ as $x\to\infty$, there exists $N_1(\epsilon)$ such that $f(x_h(n))\leq \epsilon x_h(n)$ for all $n\geq N_1(\epsilon)$. Hence for $j\geq N_1(\epsilon)$ we have
\begin{align*}
x_h(j+1)&\leq x_h(j-N)+h\sum_{n=j-N}^{j} \epsilon x_h(n) +h\sum_{n=j-N}^{j}g(x_h(n-N))\\
&\leq x_h(j-N)+ h(N+1)\epsilon x_h(j) +h\sum_{n=j-N}^{j}g(x_h(n-N))\\
&\leq x_h(j-N)+ h(N+1)\epsilon x_h(j+1) +h\sum_{n=j-N}^{j}g(x_h(n-N)).
\end{align*}
Hence for $j\geq N_1(\epsilon)$ we have
\[
x_h(j+1) \leq \frac{1}{1-(\tau+h)\epsilon}x_h(j-N)+\frac{h}{1-(\tau+h)\epsilon}h\sum_{n=j-N}^{j}g(x_h(n-N)).
\]
Since $\log g(x)/\log x \to \beta$ as $x\to\infty$, and $x_h(n)\to\infty$ as $n\to\infty$, for every $\epsilon>0$
there exists $N_2(\epsilon)\geq N$ such that $g(x_h(n-N))<x_h(n-N)^{\beta+\epsilon}$. Thus for $j\geq N_2(\epsilon)+N$ we have
\begin{align*}
h\sum_{n=j-N}^{j}g(x_h(n-N))&\leq h\sum_{n=j-N}^{j}x_h(n-N)^{\beta+\epsilon}\\
&\leq h(N+1) x_h(j-N)^{\beta+\epsilon}.
\end{align*}
Hence
\[
h\sum_{n=j-N}^{j}g(x_h(n-N)) \leq (\tau+h)x_h(j-N)^{\beta+\epsilon}, \quad j\geq N_2(\epsilon).
\]
Let $N_3(\epsilon)=\max(N_1(\epsilon),N_2(\epsilon)+N)$. Then as $ 1-(\tau+h)\epsilon>1/2$, for $j\geq N_3(\epsilon)$ we have
\[
x_h(j+1) \leq 2x_h(j-N) + 2(\tau+h) x_h(j-N)^{\beta+\epsilon}.
\]
Define $x^\ast_h(n)=x_h(n(N+1))$ for $n\geq -1$. Therefore for $n\geq N_3(\epsilon)$ we have
\begin{align*}
x^\ast_h(j+1)
&\leq x_h((j+1)(N+1))\leq 2x_h(j(N+1)) + 2(\tau+h) x_h(j(N+1))^{\beta+\epsilon}\\
&=2x^\ast_h(j) + 2(\tau+h) x^\ast_h(j)^{\beta+\epsilon}.
\end{align*}
Thus
\[
\log x^\ast_h(j+1) \leq \log 2(\tau+h) + (\beta+\epsilon)\log x^\ast_h(j)
+ \log\left( 1 + \frac{x^\ast_h(j)}{(\tau+h) x^\ast_h(j)^{\beta+\epsilon}}\right).
\]
Thus we have, with $u(n)=\log x^\ast_h(n)$, and all $n>N_5(\epsilon)$, the inequality
\[
u(n+1)\leq (\beta+2\epsilon) u(n).
\]
Thus there exists $K(\epsilon)>0$ such that $u(n)\leq K(\epsilon)(\beta+2\epsilon)^n$ for $n\geq N_5(\epsilon)$. Thus
\[
\frac{1}{n}\log u(n)\leq \frac{1}{n}\log K(\epsilon) + \log(\beta+2\epsilon).
\]
Therefore
\begin{equation*}
\limsup_{n\to\infty} \frac{\log_2 x_h(n(N+1))}{n(N+1)h}
=
\limsup_{n\to\infty} \frac{\log_2 x_h^\ast(n)}{n(N+1)h} \leq \frac{\log(\beta+2\epsilon)}{(N+1)h}=\frac{\log(\beta+2\epsilon)}{\tau+h}.
\end{equation*}
Letting $\epsilon\downarrow 0$, we arrive at
\begin{equation}   \label{eq.xNnbeta}
\limsup_{n\to\infty} \frac{\log_2 x_h(n(N+1))}{n(N+1)h} \leq \frac{\log(\beta)}{\tau+h}.
\end{equation}

For every $n\in \mathbb{N}$ there exists $j=j(n)\geq 1$ such that $n(N+1)\leq j<(n+1)(N+1)$. Since
$(x_h(n))_{n\geq 0}$ is increasing, we have
\begin{align*}
\frac{\log_2 x_h(j)}{jh}&\leq \frac{\log_2 x_h((n+1)(N+1)))}{jh} \leq \frac{\log_2 x_h((n+1)(N+1)))}{n(N+1)h}\\
&=
\frac{\log_2 x_h((n+1)(N+1)))}{(n+1)(N+1)h}\frac{n+1}{n}.
\end{align*}
By \eqref{eq.xNnbeta}, we have
\[
\limsup_{j\to\infty} \frac{\log_2 x_h(j)}{jh} \leq \frac{\log \beta}{\tau+h},
\]
which gives the desired upper limit. % in \eqref{eq.growthlimdelaydisc}.

Since $f(x)\geq 0$ we have
\[
x_h(n+1)\geq x_h(n)+hg(x_h(n-N))\geq hg(x_h(n-N))
\]
and since $x_h(n)\to\infty$ as $n\to\infty$ and $\log g(x)/\log x\to\beta$ as $x\to\infty$, it follows that for
every $\epsilon<\beta$ there exists $N_6(\epsilon)$ such that $hg(x_h(n-N))\geq x_h(n-N)^{\beta-\epsilon}>e$ for $n\geq N_5(\epsilon)$.
Hence for $n\geq N_6(\epsilon)$ we have
\[
x_h(n+1)\geq x_h(n-N))^{\beta-\epsilon}.
\]
Therefore with $u(n)=\log x_h(n)$, we have that
\[
u(n+1)=\log x_h(n+1) \geq (\beta-\epsilon) u(n-N).
\]
Therefore, there exists $k(\epsilon)>0$ such that $u(n)\geq k(\epsilon) (\beta-\epsilon)^{n/(N+1)}$ for $n\geq N_6(\epsilon)$. Therefore
\[
\frac{1}{n}\log u(n)\geq \frac{1}{n}\log k(\epsilon) + \frac{1}{N+1}\log (\beta-\epsilon).
\]
Hence
\[
\liminf_{n\to\infty} \frac{\log_2 x_h(n)}{nh} \geq \frac{\log(\beta-\epsilon)}{(N+1)h}=\frac{\log(\beta-\epsilon)}{\tau+h}.
\]
Letting $\epsilon\downarrow 0$, we get
\[
\liminf_{n\to\infty} \frac{\log_2 x_h(n)}{nh} \geq \frac{\log \beta}{\tau+h}.
\]
and so combining this with the other limit we get
\[
\lim_{n\to\infty} \frac{\log_2 x_h(n)}{nh} =\frac{\log\beta}{\tau+h},
\]
as required.

The proof that \eqref{eq.growthlimdelaydiscrvbetainterp} follows from \eqref{eq.growthlimdelaydiscrvbeta} is identical in all regards to the proof of Theorem~\ref{thm:presrv1} that \eqref{eq.growthlimdelaydiscrv1interp} follows from
\eqref{eq.growthlimdelaydisc}, and is therefore omitted.

\end{document}